\documentclass[12pt,a4paper]{article}

\usepackage[margin=2.5cm]{geometry}
\usepackage[T1]{fontenc}
\usepackage[utf8]{inputenc}
\usepackage{lmodern}
\usepackage{microtype}
\usepackage{amsmath,amssymb,amsthm}
\usepackage{indentfirst}
\usepackage{xcolor}
\usepackage{cite}
\usepackage{needspace}
\usepackage[bottom]{footmisc}
\usepackage[
colorlinks=true,
linkcolor=black,
citecolor=black,
urlcolor=blue
]{hyperref}

\numberwithin{equation}{section}

\theoremstyle{plain}
\newtheorem{theorem}{Theorem}[section]
\newtheorem{proposition}[theorem]{Proposition}

\newtheorem{corollary}[theorem]{Corollary}

\theoremstyle{definition}
\newtheorem{remark}[theorem]{Remark}

\newcommand{\R}{\mathbb R}

\hypersetup{
  pdftitle={Stability for the Navier-Stokes and Incompressible Porous Media Equations with Only Vertical Dissipation on R2},
  pdfauthor={Marius Paicu and Ning Zhu}
}
\date{}

\begin{document}

\title{Stability for the Navier--Stokes and
Incompressible Porous Media Equations with Only Vertical Dissipation
on $\R^2$}
\author{
	Marius Paicu\footnote{Institut de Math\'ematiques de Bordeaux, Universit\'e de Bordeaux, F-33405 Talence Cedex, France. Email: marius.paicu@math.u-bordeaux.fr}
	\and
	Ning Zhu\footnote{School of Mathematics, Shandong University, Jinan, Shandong 250100, P. R. China. \\Email: ning.zhu@sdu.edu.cn}
}
\maketitle

\begin{abstract}
For the two-dimensional Navier--Stokes equations on $\R^2$, we prove that
vertical viscosity alone yields global $H^2$ stability and algebraic decay
of the vorticity and velocity for small reflection-symmetric data with
polynomially weighted vorticity.
Under stronger smallness assumptions, we obtain $H^3$ stability and
vorticity-gradient decay,
in contrast to unbounded gradient growth for Euler flows in the same
symmetry class.
For incompressible porous media (IPM), we obtain $H^2$ stability and decay
under analogous assumptions with only vertical diffusion, despite the
absence of derivative gain in Darcy's law.
This complements finite-time blow-up constructions for nondiffusive IPM
in periodic domains.
The proofs combine a weighted inequality for odd functions with anisotropic
energy estimates exploiting the structure of each system.
\end{abstract}



\section{Introduction}

The long-time behavior of two-dimensional incompressible flows depends
strongly on the available dissipation. Smooth Euler solutions can develop
unbounded vorticity gradients, with double-exponential growth in domains
with a boundary \cite{KS,ZlatosHalfPlane} and linear or superlinear growth
examples in the whole plane \cite{CJ,JYaoZhou}. Full Laplacian viscosity,
on the other hand, yields decay of finite-energy solutions and, under
suitable integrability or weighted assumptions, of higher-order norms
\cite{Schonbek,Wiegner,SCho4,SW}. A natural question is whether viscosity
in a single direction can stabilize small perturbations of the zero
solution and prevent growth of higher Sobolev norms.

We first consider the two-dimensional incompressible Navier--Stokes
equations with only vertical viscosity,
\begin{equation}\label{eq:velocity-system}
\left\{
\begin{aligned}
&\partial_tu+u\cdot\nabla u-\nu\partial_y^2u=-\nabla p,
\qquad (x,y)\in\R^2,\quad t>0,\\
&\partial_xu^1+\partial_yu^2=0,\\
&u|_{t=0}=u_0.
\end{aligned}
\right.
\end{equation}
Here $u=(u^1,u^2)$ is the velocity, $p$ is the pressure, and $\nu>0$.
We take $\nu=1$ below. The general case follows from
$u(t,x,y)=\nu v(\nu t,x,y)$ and $p(t,x,y)=\nu^2q(\nu t,x,y)$,
where $(v,q)$ solves the normalized system.

Directional diffusion is relevant to flows with widely separated
horizontal and vertical length scales. Under an anisotropic rescaling,
the coefficients of the two diffusive terms differ by the square of
the aspect ratio. A small vertical scale therefore leads to dominant
vertical diffusion. The Prandtl boundary-layer equations provide a
classical example of this balance: normal diffusion is retained in the
tangential momentum equation, while tangential diffusion is neglected
\cite{Pran,Sch}. System \eqref{eq:velocity-system} shares this directional
viscous structure, while retaining both momentum equations and the full
pressure coupling. It isolates a basic question in anisotropic fluid
models: can diffusion in one direction control the nonlinear creation
of derivatives in the other?

In three dimensions, the equations with horizontal viscosity and no
vertical viscosity have a well-developed small-data global theory
\cite{CDGG,Paicu,CZ,PZ}. Under additional low-frequency or integrability
assumptions, stability and decay estimates have also been established
\cite{JWY,XZ,Fujii}. Both these equations and
\eqref{eq:velocity-system} leave one direction undamped. However, the
three-dimensional problem retains diffusion on a two-dimensional plane,
whereas the planar problem has only one dissipative variable. This
reduction weakens both the interpolation estimates for transport and
the linear decay.

For the two-dimensional equations, Liang, Zhang and Zhu \cite{LZZ}
established global well-posedness in an anisotropic Sobolev space.
Zhou and Wu \cite{ZW} proved global existence and uniqueness in
$H^1(\R^2)$, together with the corresponding energy identity.
At the $H^2$ level, our analysis couples uniform control of horizontal
vorticity derivatives to quantitative decay from vertical viscosity.
On
$\mathbb T\times\R$, Dong, Wu, Xu and Zhu \cite{Dong} proved $H^2$
stability with horizontal dissipation and exponential decay of the
oscillatory velocity in $H^1$. Cao and Guo \cite{CG} studied
well-posedness and long-time behavior in a horizontally periodic
channel. Periodicity in the dissipative direction provides a
Poincar\'e inequality for the oscillatory part, with the mean mode
treated separately. On $\R^2$, arbitrarily low frequencies remain in
the dissipative direction, and no such uniform bound is available.
For fractional horizontal dissipation of order $2s$, Wang,
Wu and Zhu \cite{WangWuZhu} obtained small-data stability on $\R^2$
for $0\leq s<1$ under assumptions depending on $s$. The ordinary
second-order case considered here lies outside that range.

The vorticity formulation separates the transport structure from the
pressure and relates the velocity derivatives to a scalar quantity. We
use the convention
\begin{equation*}
\omega=\partial_yu^1-\partial_xu^2.
\end{equation*}
For a divergence free velocity, the $L^2$ norm of $\nabla u$ equals
that of $\omega$. The vorticity satisfies the transport--diffusion equation
\begin{equation}\label{eq:vorticity-equation}
\partial_t\omega+u\cdot\nabla\omega-\partial_y^2\omega=0.
\end{equation}

The main difficulty in obtaining a uniform $H^2$ bound is already
visible in the energy identity for one vorticity derivative:
\begin{equation*}\begin{aligned}
&\frac12\frac{d}{dt}\|\nabla\omega\|_{L^2}^2
+\|\partial_y\nabla\omega\|_{L^2}^2\\
={} &-\int_{\R^2}\partial_xu^1(\partial_x\omega)^2\,dx\,dy
-\int_{\R^2}\partial_yu^1\partial_x\omega\partial_y\omega\,dx\,dy\\
&\quad-\int_{\R^2}\partial_xu^2\partial_y\omega\partial_x\omega\,dx\,dy
-\int_{\R^2}\partial_yu^2(\partial_y\omega)^2\,dx\,dy.
\end{aligned}\end{equation*}
The first integral contains two purely horizontal vorticity
derivatives. Incompressibility, anisotropic H\"older's inequality
and interpolation give the bound
\[
\left|\int_{\R^2}\partial_xu^1(\partial_x\omega)^2\,dx\,dy\right|
\leq\frac18\|\partial_{xy}\omega\|_{L^2}^2
+C\|u\|_{L^2}^2\|\omega\|_{L^2}^2\|\partial_x\omega\|_{L^2}^2.
\]
The first term is absorbed by the vertical dissipation, whereas the
remaining horizontal derivative is not controlled by this dissipation
on $\R^2$. Together with the other nonlinear terms, this leads to
\[
\frac{d}{dt}\|\nabla\omega\|_{L^2}^2
+\|\partial_y\nabla\omega\|_{L^2}^2
\leq C\bigl(1+\|u_0\|_{L^2}^2\bigr)
\|\omega\|_{L^2}^2\|\nabla\omega\|_{L^2}^2.
\]
Time integrability of the squared vorticity norm therefore suffices
to close this estimate. The basic vorticity energy identity controls
its vertical derivative in time, but does not provide this
integrability of the vorticity itself.

The linear decay in one dissipative variable does not close this
estimate either. For data in $L_x^2L_y^1$, without a cancellation
condition, the heat estimate is
\[
\|e^{t\partial_y^2}f\|_{L^2}
\leq Ct^{-\frac14}\|f\|_{L_x^2L_y^1},\qquad t>0.
\]
The square of this bound is not integrable at infinity. Thus the
usual linear estimate, without additional cancellation, does not
provide the integrable coefficient needed in this energy argument.

We obtain the required decay by combining the preserved reflection
symmetry with a polynomially weighted $L^2$ condition on the initial
vorticity in the vertical variable. Oddness of the vorticity yields a
lower bound for the weighted energy at the vertical heat scale. The
resulting decay exponent exceeds $1/2$, which closes the $H^2$
estimate. The same weighted bound controls a vertical moment and
also yields decay of the velocity.

For $m>2$, define
\begin{equation}\label{eq:weight}
w_m(y)=\left(1+\frac{y^2}{2(2m-1)}\right)^{\frac m2}.
\end{equation}
Our Navier--Stokes result is the following.

\begin{theorem}[Stability and decay for Navier--Stokes]
\label{thm:main}
Let $m>2$. Suppose that $u_0\in H^2(\R^2)$ is divergence free and
\begin{equation}\label{eq:initial-symmetry}
u_0^1(x,-y)=u_0^1(x,y),\qquad
u_0^2(x,-y)=-u_0^2(x,y).
\end{equation}
There exists $\varepsilon_m>0$ such that, if
\begin{equation*}
\|u_0\|_{H^2}+\|w_m\omega_0\|_{L^2}
\leq\varepsilon\leq\varepsilon_m,
\end{equation*}
then \eqref{eq:velocity-system} has a unique global
solution in $C([0,\infty);H^2(\R^2))$ satisfying
\begin{equation*}
\sup_{t\geq0}\|u(t)\|_{H^2}^2
+\int_0^\infty\|\partial_yu(t)\|_{H^2}^2\,dt
\leq C_m\varepsilon^2.
\end{equation*}
For every $t\geq0$,
\begin{align}
\left\|w_m\!\left(\frac{y}{\sqrt{1+t}}\right)\omega(t)\right\|_{L^2}
&\leq C_m\|w_m\omega_0\|_{L^2}
(1+t)^{-\frac{3m}{4(m+1)}},
\label{eq:main-decay-estimate}\\*
\|u(t)\|_{L^2}
&\leq C_m\|w_m\omega_0\|_{L^2}
(1+t)^{-\frac{m-2}{4(m+1)}}.
\label{eq:main-velocity-decay}
\end{align}
In particular, \eqref{eq:main-decay-estimate} holds with
$\|\omega(t)\|_{L^2}$ on the left-hand side.

If, in addition, $u_0\in H^3(\R^2)$ and the stronger smallness
condition
\begin{equation}\label{eq:H3-small-data}
\|u_0\|_{H^3}+\|w_m\omega_0\|_{L^2}
\leq\varepsilon\leq\varepsilon_m
\end{equation}
holds, then the solution belongs to $C([0,\infty);H^3(\R^2))$ and
\begin{equation}\label{eq:main-H3-estimate}
\sup_{t\geq0}\|u(t)\|_{H^3}^2
+\int_0^\infty\|\partial_yu(t)\|_{H^3}^2\,dt
\leq C_m\varepsilon^2.
\end{equation}
\end{theorem}

\Needspace{7\baselineskip}
Reflection across the boundary gives the corresponding statement
with impermeable slip boundary conditions on the upper half-plane.

\begin{corollary}[The upper half-plane]\label{cor:half-plane}
Let $m>2$ and $\R^2_+=\R\times(0,\infty)$. Consider
\eqref{eq:velocity-system} on $\R^2_+$ with
\begin{equation}\label{eq:half-plane-boundary}
u^2(t,x,0)=0,\qquad \partial_yu^1(t,x,0)=0.
\end{equation}
Suppose that $u_0\in H^2(\R^2_+)$ is divergence free, satisfies these
boundary conditions, and
\[
\|u_0\|_{H^2(\R^2_+)}+\|w_m\omega_0\|_{L^2(\R^2_+)}
\leq\varepsilon\leq\varepsilon_m.
\]
Then the conclusions of Theorem~\ref{thm:main} hold with all norms
taken over $\R^2_+$. If the same assumptions hold with $H^3$ in
place of $H^2$, then \eqref{eq:main-H3-estimate} also holds on
$\R^2_+$.
\end{corollary}

\begin{remark}\label{rem:gaussian}
For Navier--Stokes, replacing the polynomial weighted assumption by
smallness of $\|u_0\|_{H^2}+\|e^{y^2/8}\omega_0\|_{L^2}$ improves
the vorticity and velocity decay rates to $(1+t)^{-3/4}$ and
$(1+t)^{-1/4}$, respectively. This follows from the Gaussian-weight
calculation at the end of Subsection~\ref{subsec:proof-coercivity}.
\end{remark}

\begin{remark}
If $\|u(t)\|_{\dot H^{-\sigma}}$ is uniformly bounded for some
$\sigma>0$, interpolation gives
\[
\|u(t)\|_{L^2}
\leq\|u(t)\|_{\dot H^{-\sigma}}^{\frac1{1+\sigma}}
\|\omega(t)\|_{L^2}^{\frac{\sigma}{1+\sigma}}.
\]
This improves \eqref{eq:main-velocity-decay} when
$\sigma>(m-2)/[2(m+1)]$, provided that the negative Sobolev norm
can be propagated uniformly in time.
\end{remark}

\begin{remark}[Comparison with inviscid growth]\label{rem:euler}
Choi and Jeong \cite[Theorem~1.2]{CJ} constructed smooth, compactly
supported initial vorticities, odd in the vertical variable, whose
Euler solutions have unbounded $L^2$ vorticity gradients. Amplitude
and time scaling makes their initial velocities arbitrarily small in
the norms in \eqref{eq:H3-small-data}. For the vertically dissipative
solution starting from the same data, our estimates and interpolation
give
\[
\|u(t)\|_{\dot H^2}
=\|\nabla\omega(t)\|_{L^2}
\leq\|\omega(t)\|_{L^2}^{\frac{1}{2}}
\|\omega(t)\|_{\dot H^2}^{\frac{1}{2}}
\leq C_m\varepsilon(1+t)^{-\frac{3m}{8(m+1)}}.
\]
Thus the reflection symmetry is compatible with inviscid growth,
whereas vertical viscosity yields decay for these small data.
\end{remark}

The same mechanism applies to the incompressible porous media (IPM) equation
with vertical diffusion,
\begin{equation}\label{eq:IPM-system}
\left\{
\begin{aligned}
&\partial_t\rho+u\cdot\nabla\rho-\partial_y^2\rho=0,\\
&u=-\nabla p-\begin{pmatrix}0\\\rho\end{pmatrix},\\
&\partial_xu^1+\partial_yu^2=0,\\
&\rho|_{t=0}=\rho_0,
\end{aligned}
\right.
\qquad (x,y)\in\R^2,\quad t>0.
\end{equation}
Here $\rho$ is the density fluctuation about a constant state, and $u$
is determined by Darcy's law. In porous-medium transport, anisotropic
thermal or solute diffusion leads to directional diffusivities.
Equation \eqref{eq:IPM-system} represents the limiting case in which
diffusion acts in one direction only. Here we obtain uniform $H^2$ stability and algebraic decay for the
planar equation with a single dissipative variable, under reflection
symmetry and a polynomially weighted density condition.

For the nondiffusive equation, Kiselev and Yao \cite{KY} constructed
solutions whose derivatives must grow without bound if they remain
globally smooth. C\'ordoba and Mart\'inez-Zoroa \cite{CMIPM}
constructed finite-time singularities for smooth solutions with a
smooth source. Recent preprints of Dai \cite{DaiIPM} and Zheng
\cite{ZhengIPM} give finite-time blow-up constructions for the
unforced equation on $\mathbb T^2$.
These results further motivate the stabilization question addressed
here for the vertically diffusive equation on $\R^2$.

In \eqref{eq:IPM-system}, the velocity is an operator
of order zero on $\rho$. Thus the derivative gained by the
Biot--Savart law in \eqref{eq:vorticity-equation} is unavailable.
We use directional identities from Darcy's law to recover vertical
derivatives in the highest-order horizontal interactions.
For the weighted transport term, the decomposition
$u^2=-\rho+(u^2+\rho)$ isolates a remainder whose $L^\infty$ norm
can be controlled using only vertical derivatives of $\rho$.
These two features allow the density energy and weighted decay
estimates to close together at the $H^2$ level. The two results thus
address the same stabilization question for velocity laws of
different orders.

With the same polynomial weight \eqref{eq:weight}, we obtain the
following stability and decay result for the density.

\begin{theorem}[$H^2$ stability and decay for IPM]\label{thm:IPM}
Let $m>2$. Suppose that $\rho_0\in H^2(\R^2)$ is odd in $y$ and
\begin{equation*}
\|\rho_0\|_{H^2}+\|w_m\rho_0\|_{L^2}
\leq\varepsilon\leq\varepsilon_m.
\end{equation*}
For $\varepsilon_m>0$ sufficiently small, \eqref{eq:IPM-system} has
a unique global solution $\rho\in C([0,\infty);H^2(\R^2))$ such that
\begin{equation}\label{eq:IPM-main-energy}
\sup_{t\geq0}\|\rho(t)\|_{H^2}^2
+\int_0^\infty\|\partial_y\rho(t)\|_{H^2}^2\,dt
\leq C_m\varepsilon^2.
\end{equation}
Moreover,
\begin{align}
\left\|w_m\!\left(\frac{y}{\sqrt{1+t}}\right)\rho(t)\right\|_{L^2}
&\leq C_m\|w_m\rho_0\|_{L^2}
(1+t)^{-\frac{3m}{4(m+1)}},
\label{eq:IPM-main-decay}\\*
\|u(t)\|_{L^2}
&\leq C_m\|w_m\rho_0\|_{L^2}
(1+t)^{-\frac{3m}{4(m+1)}}.
\label{eq:IPM-main-velocity}
\end{align}
\end{theorem}

The proofs use weighted scalar estimates together with anisotropic
energy estimates. Gallay and Wayne \cite{GW} developed polynomially
weighted vorticity spaces and self-similar variables to study the
long-time behavior of the fully viscous Navier--Stokes equations.
Here the diffusion acts only in $y$: the weighted estimate provides
no horizontal smoothing, while the differentiated transport contains
purely horizontal derivatives. We therefore combine a one-dimensional
lower bound on odd functions with estimates that control the
horizontal interactions through vertical dissipation and time decay.

The rest of the paper is organized as follows.
Section~\ref{sec:NS} proves the Navier--Stokes theorem, starting with
the $H^2$ energy estimate and then establishing weighted vorticity
decay, velocity decay, and the higher-order estimate. The half-plane
corollary is proved at the end of that section.
Section~\ref{sec:IPM} treats IPM, using the common weighted identity
and closing the density energy and decay estimates simultaneously.

Throughout the paper, spatial norms are taken over $\R^2$ unless
another domain is specified. In the mixed norm $L_x^pL_y^q$, the
$L_y^q$ norm is taken first. The notation $\approx$ denotes two-sided
bounds with constants independent of time. The constants $C$ and
$C_m$ may change from line to line; $C_m$ may depend on $m$.

\section{Stability and decay for Navier--Stokes}\label{sec:NS}

We first reduce the uniform $H^2$ bound to the time-integrability of
the squared vorticity norm. The weighted estimate then supplies this
integrability and yields velocity decay. Once these bounds are
available, only the homogeneous third-order velocity estimate remains
to prove the higher-regularity assertion.

\subsection{The \texorpdfstring{$H^2$}{H2} energy estimate and bootstrap}
\label{subsec:H2-bootstrap}

Let $T>0$ be contained in the interval of smooth existence, and define
\begin{equation*}
	E_2(T)
	=
	\sup_{0\leq t\leq T}\|u(t)\|_{H^2}^2
	+\int_0^T\|\partial_yu(t)\|_{H^2}^2\,dt.
\end{equation*}
Let $C_0\geq1$ be a fixed constant, specified below independently of
$\varepsilon$ and $T$. Assume that
\begin{equation}\label{eq:E2-bootstrap}
	E_2(T)\leq4C_0\varepsilon^2.
\end{equation}

To close this bootstrap, it is enough to establish the energy inequality
\begin{align}
	E_2(T)
	&\leq C\|u_0\|_{H^2}^2
	+C\bigl(1+\|u_0\|_{L^2}^2\bigr)
	\int_0^T
	\|\omega(t)\|_{L^2}^2
	\|u(t)\|_{H^2}^2\,dt
\label{eq:H2-energy-reduction}
\end{align}
and the decay estimate
\begin{equation}\label{eq:decay-reduction}
	\|\omega(t)\|_{L^2}
	\leq
	C_m\|w_m\omega_0\|_{L^2}
	(1+t)^{-\frac{3m}{4(m+1)}},
	\qquad 0\leq t\leq T.
\end{equation}
Indeed, since $m>2$, \eqref{eq:decay-reduction} implies
\begin{equation*}
	\int_0^T\|\omega(t)\|_{L^2}^2\,dt
	\leq
	C_m\|w_m\omega_0\|_{L^2}^2
	\leq C_m\varepsilon^2.
\end{equation*}
Applying Gronwall's inequality to \eqref{eq:H2-energy-reduction}, we obtain
\[
E_2(T)
\leq
C\|u_0\|_{H^2}^2
\exp\!\left\{
C_m\bigl(1+\|u_0\|_{L^2}^2\bigr)\varepsilon^2
\right\}.
\]
Choose $C_0$ to be twice the constant multiplying the initial energy
in this bound, increasing it to $1$ if necessary. Then choose
$\varepsilon_m$ so that the exponential is at most $2$. It follows that
\[
E_2(T)\leq2C_0\varepsilon^2,
\]
which improves \eqref{eq:E2-bootstrap}. The standard continuity argument
therefore closes the bootstrap and proves the global $H^2$ estimate.

We establish \eqref{eq:H2-energy-reduction} first; the next
subsection proves \eqref{eq:decay-reduction}.
Taking the $L^2$ inner product of \eqref{eq:velocity-system} with $u$, we
obtain
\begin{equation}\label{eq:u-L2}
	\|u(t)\|_{L^2}^2
	+2\int_0^t\|\partial_yu(s)\|_{L^2}^2\,ds
	=\|u_0\|_{L^2}^2.
\end{equation}
Similarly, taking the $L^2$ inner product of
\eqref{eq:vorticity-equation} with $\omega$ gives
\begin{equation}\label{eq:omega-L2}
	\|\omega(t)\|_{L^2}^2
	+2\int_0^t\|\partial_y\omega(s)\|_{L^2}^2\,ds
	=\|\omega_0\|_{L^2}^2.
\end{equation}

Differentiating \eqref{eq:vorticity-equation} once and taking the
corresponding $L^2$ inner products, we find
\begin{align}
	&\frac12\frac{d}{dt}\|\nabla\omega\|_{L^2}^2
	+\|\partial_y\nabla\omega\|_{L^2}^2\notag\\*
	={} &
	-\int_{\R^2}\partial_xu^1(\partial_x\omega)^2\,dx\,dy
	-\int_{\R^2}\partial_yu^1\partial_x\omega\partial_y\omega\,dx\,dy\notag\\*
	&\quad-\int_{\R^2}\partial_xu^2\partial_y\omega\partial_x\omega\,dx\,dy
	-\int_{\R^2}\partial_yu^2(\partial_y\omega)^2\,dx\,dy\notag\\*
	=:{} &A_1+A_2+A_3+A_4.
\label{eq:gradient-vorticity}
\end{align}
We first estimate $A_1$. Since
$\partial_xu^1=-\partial_yu^2$, integration by parts in $y$ gives
\[
A_1
=-2\int_{\R^2}u^2\partial_{xy}\omega\partial_x\omega\,dx\,dy.
\]
Using anisotropic H\"older's inequality, one-dimensional interpolation,
and incompressibility, we obtain
\begin{align}
	|A_1|
	&\leq
	C\|u^2\|_{L_y^2L_x^\infty}
	\|\partial_x\omega\|_{L_x^2L_y^\infty}
	\|\partial_{xy}\omega\|_{L^2}\notag\\
	&\leq
	C\|u^2\|_{L^2}^{\frac{1}{2}}
	\|\partial_xu^2\|_{L^2}^{\frac{1}{2}}
	\|\partial_x\omega\|_{L^2}^{\frac{1}{2}}
	\|\partial_{xy}\omega\|_{L^2}^{\frac{3}{2}}\notag\\
	&\leq
	C\|u\|_{L^2}^{\frac{1}{2}}
	\|\omega\|_{L^2}^{\frac{1}{2}}
	\|\nabla\omega\|_{L^2}^{\frac{1}{2}}
	\|\partial_y\nabla\omega\|_{L^2}^{\frac{3}{2}}\notag\\
	&\leq
	\frac18\|\partial_y\nabla\omega\|_{L^2}^2
	+C\|u\|_{L^2}^2\|\omega\|_{L^2}^2
	\|\nabla\omega\|_{L^2}^2.
	\label{eq:A1-estimate}
\end{align}

For $A_2$, anisotropic H\"older's inequality in the two spatial directions
gives
\begin{align}
	|A_2|
	&\leq
	\|\partial_yu^1\|_{L^2}
	\|\partial_x\omega\|_{L_x^2L_y^\infty}
	\|\partial_y\omega\|_{L_y^2L_x^\infty}\notag\\*
	&\leq
	C\|\omega\|_{L^2}
	\|\partial_x\omega\|_{L^2}^{\frac{1}{2}}
	\|\partial_{xy}\omega\|_{L^2}^{\frac{1}{2}}
	\|\partial_y\omega\|_{L^2}^{\frac{1}{2}}
	\|\partial_{xy}\omega\|_{L^2}^{\frac{1}{2}}\notag\\*
	&\leq
	\frac18\|\partial_y\nabla\omega\|_{L^2}^2
	+C\|\omega\|_{L^2}^2\|\nabla\omega\|_{L^2}^2.
	\label{eq:A2-estimate}
\end{align}
The term $A_3$ is estimated in exactly the same way:
\begin{equation}\label{eq:A3-estimate}
	|A_3|
	\leq
	\frac18\|\partial_y\nabla\omega\|_{L^2}^2
	+C\|\omega\|_{L^2}^2\|\nabla\omega\|_{L^2}^2.
\end{equation}

Finally, using $\partial_yu^2=-\partial_xu^1$ and integrating by parts in
$x$, we have
\[
A_4
=-2\int_{\R^2}u^1\partial_{xy}\omega\partial_y\omega\,dx\,dy.
\]
Arguing as in the estimate of $A_1$, with $u^1$ in
$L_x^2L_y^\infty$ and $\partial_y\omega$ in
$L_y^2L_x^\infty$, we obtain
\begin{align}
	|A_4|
	&\leq
	C\|u^1\|_{L_x^2L_y^\infty}
	\|\partial_y\omega\|_{L_y^2L_x^\infty}
	\|\partial_{xy}\omega\|_{L^2}\notag\\*
	&\leq C\|u^1\|_{L^2}^{\frac{1}{2}}
	\|\partial_yu^1\|_{L^2}^{\frac{1}{2}}
	\|\partial_y\omega\|_{L^2}^{\frac{1}{2}}
	\|\partial_{xy}\omega\|_{L^2}^{\frac{3}{2}}\notag\\*
	&\leq
	\frac18\|\partial_y\nabla\omega\|_{L^2}^2
	+C\|u\|_{L^2}^2\|\omega\|_{L^2}^2
	\|\nabla\omega\|_{L^2}^2.
	\label{eq:A4-estimate}
\end{align}

Combining \eqref{eq:gradient-vorticity} with
\eqref{eq:A1-estimate}--\eqref{eq:A4-estimate} and using
\eqref{eq:u-L2}, we arrive at
\begin{equation}\label{eq:H2-differential}
	\frac{d}{dt}\|\nabla\omega\|_{L^2}^2
	+\|\partial_y\nabla\omega\|_{L^2}^2
	\leq
	C\bigl(1+\|u_0\|_{L^2}^2\bigr)
	\|\omega\|_{L^2}^2\|\nabla\omega\|_{L^2}^2.
\end{equation}

The Biot--Savart law gives
\begin{equation*}
	\|\nabla\partial_x^a\partial_y^b u\|_{L^2}
	=\|\partial_x^a\partial_y^b\omega\|_{L^2},
\qquad a,b\in\{0,1,2,\ldots\}.
\end{equation*}

In particular,
\begin{equation}\label{eq:H2-equivalence}
	\|u\|_{H^2}^2
	\approx
	\|u\|_{L^2}^2+\|\omega\|_{H^1}^2,
	\qquad
	\|\partial_yu\|_{H^2}^2
	\approx
	\|\partial_yu\|_{L^2}^2
	+\|\partial_y\omega\|_{H^1}^2.
\end{equation}
Integrating \eqref{eq:H2-differential} over $[0,T]$ and combining it with
\eqref{eq:u-L2}, \eqref{eq:omega-L2}, and
\eqref{eq:H2-equivalence}, we obtain
\begin{align*}
	E_2(T)
	&\leq
	C\|u_0\|_{H^2}^2
	+C\bigl(1+\|u_0\|_{L^2}^2\bigr)
	\int_0^T
	\|\omega(t)\|_{L^2}^2
	\|u(t)\|_{H^2}^2\,dt.
\end{align*}
This is precisely \eqref{eq:H2-energy-reduction}.

\subsection{Weighted vorticity decay}\label{sec:decay}

To obtain the decay required in \eqref{eq:decay-reduction}, we use a
polynomial weight that expands at the vertical heat scale. We first
record the weighted equation for a transported scalar, so that the
same identity can later be applied to the IPM density.
Consider
\begin{equation}\label{eq:scalar-equation}
\partial_tf+v\cdot\nabla f-\partial_y^2f=0,
\qquad \nabla\cdot v=0.
\end{equation}
The self-similar vertical variable is $\eta=y/\sqrt{1+t}$.
We incorporate this change of scale directly into the weight and set
\begin{equation*}
W(t,x,y)=w_m\!\left(\frac{y}{\sqrt{1+t}}\right)f(t,x,y).
\end{equation*}
The derivatives of the weight satisfy
\begin{align*}
\frac{\partial_y[w_m(\frac{y}{\sqrt{1+t}})]}{w_m(\frac{y}{\sqrt{1+t}})}
&=\frac{my}{2(2m-1)(1+t)+y^2},
\\*
\frac{\partial_t[w_m(\frac{y}{\sqrt{1+t}})]}{w_m(\frac{y}{\sqrt{1+t}})}
&=-\frac{my^2}{2(1+t)[2(2m-1)(1+t)+y^2]}.
\end{align*}
Substitution into \eqref{eq:scalar-equation} gives
\begin{equation}\begin{aligned}
&\partial_tW+v\cdot\nabla W-\partial_y^2W
+\frac{2my}{2(2m-1)(1+t)+y^2}\partial_yW
\\
&\quad+\Biggl\{
\partial_y\!\left(\frac{my}{2(2m-1)(1+t)+y^2}\right)
-\frac{m^2y^2}{[2(2m-1)(1+t)+y^2]^2}
\\
&\qquad
+\frac{my^2}{2(1+t)[2(2m-1)(1+t)+y^2]}
-\frac{myv^2}{2(2m-1)(1+t)+y^2}
\Biggr\}W=0.
\end{aligned}\label{eq:moving-weight-equation}\end{equation}
Taking the $L^2$ inner product with $W$, we use
\[
\int_{\R^2}\frac{2my}{2(2m-1)(1+t)+y^2}\partial_yW\,W\,dx\,dy
=-\int_{\R^2}\partial_y\!\left(\frac{my}{2(2m-1)(1+t)+y^2}\right)
|W|^2\,dx\,dy.
\]
This cancels the corresponding zeroth-order term in
\eqref{eq:moving-weight-equation}. Since $\nabla\cdot v=0$, we obtain
\begin{align}
&\frac12\frac{d}{dt}\|W\|_{L^2}^2
+\|\partial_yW\|_{L^2}^2
+\frac1{1+t}\int_{\R^2}
V_m\!\left(\frac{y}{\sqrt{1+t}}\right)|W|^2\,dx\,dy\notag\\*
={} &\int_{\R^2}
\frac{myv^2}{2(2m-1)(1+t)+y^2}|W|^2\,dx\,dy,
\label{eq:common-weighted-energy}
\end{align}
where
\begin{equation}\label{eq:potential}
V_m(\eta)
=\frac12\eta\frac{w_m'(\eta)}{w_m(\eta)}
-\left(\frac{w_m'(\eta)}{w_m(\eta)}\right)^2
=\frac{m\eta^2[2(m-1)+\eta^2]}
{2[2(2m-1)+\eta^2]^2}\geq0.
\end{equation}

The next step is to bound the two positive terms on the left-hand
side of \eqref{eq:common-weighted-energy} by the weighted $L^2$
norm. The reflection symmetry makes the scalar odd in $y$; the
following inequality gives the required lower bound in that class.

\begin{proposition}\label{prop:coercive-inequality}
Let $m>2$. Every odd $f\in H^1(\R)$ satisfies
\begin{equation*}
\int_\R\bigl(|f'(\eta)|^2+V_m(\eta)|f(\eta)|^2\bigr)\,d\eta
\geq\frac{3m}{4(m+1)}\int_\R|f(\eta)|^2\,d\eta.
\end{equation*}
\end{proposition}

We prove this proposition in Subsection~\ref{subsec:proof-coercivity}.
For now, apply it at each fixed $x$, with the change of variables
$\eta=y/\sqrt{1+t}$. For $W$ odd in $y$, integration in $x$ gives
\begin{equation}\label{eq:moving-coercivity}
\|\partial_yW\|_{L^2}^2
+\frac1{1+t}\int_{\R^2}
V_m\!\left(\frac{y}{\sqrt{1+t}}\right)|W|^2\,dx\,dy
\geq\frac{3m}{4(m+1)(1+t)}\|W\|_{L^2}^2.
\end{equation}
The coefficient $3m/[4(m+1)]$ exceeds $1/2$ for $m>2$.
It remains to control the transport contribution while retaining
this coefficient in the decay estimate.

We apply \eqref{eq:common-weighted-energy} with $f=\omega$ and $v=u$.
The reflection
\[
(u^1,u^2)(t,x,y)\mapsto(u^1,-u^2)(t,x,-y)
\]
preserves \eqref{eq:velocity-system}. Uniqueness and
\eqref{eq:initial-symmetry} therefore imply that $u^1$ remains even
and $u^2$ remains odd in $y$. Consequently, $\omega$ and the weighted
vorticity $W$ are odd, so \eqref{eq:moving-coercivity} applies.

From \eqref{eq:potential},
\[
\left(\frac{m\eta}{2(2m-1)+\eta^2}\right)^2
\leq\frac{m}{m-1}V_m(\eta).
\]
Anisotropic H\"older's inequality and one-dimensional interpolation
give
\begin{align}
&\left|\int_{\R^2}
\frac{myu^2}{2(2m-1)(1+t)+y^2}|W|^2\,dx\,dy\right|\notag\\
\leq{} &
\left\|\frac{myW}{2(2m-1)(1+t)+y^2}\right\|_{L^2}
\|u^2\|_{L_y^2L_x^\infty}\|W\|_{L_x^2L_y^\infty}\notag\\
\leq{} & C_m\left(\frac1{1+t}\int_{\R^2}
V_m\!\left(\frac{y}{\sqrt{1+t}}\right)|W|^2\,dx\,dy\right)^{\frac12}
\|u^2\|_{L^2}^{\frac12}\|\partial_xu^2\|_{L^2}^{\frac12}
\|W\|_{L^2}^{\frac12}\|\partial_yW\|_{L^2}^{\frac12}.
\label{eq:NS-weighted-holder}
\end{align}
The velocity energy identity, the Biot--Savart law, and $w_m\geq1$
imply
\[
\|u^2(t)\|_{L^2}\leq\|u_0\|_{L^2},
\qquad
\|\partial_xu^2(t)\|_{L^2}\leq\|\omega(t)\|_{L^2}
\leq\|W(t)\|_{L^2}.
\]
Thus the two undifferentiated factors in
\eqref{eq:NS-weighted-holder} contribute one power of
$\|W\|_{L^2}$. The potential term contributes one half-power of
the quadratic form, and the derivative factor contributes one
quarter-power. Consequently,
\begin{align*}
&\left|\int_{\R^2}
\frac{myu^2}{2(2m-1)(1+t)+y^2}|W|^2\,dx\,dy\right|\\
\leq{} & C_m\|u_0\|_{L^2}^{\frac12}\|W\|_{L^2}
\left(\frac1{1+t}\int_{\R^2}
V_m\!\left(\frac{y}{\sqrt{1+t}}\right)|W|^2\,dx\,dy\right)^{\frac12}
\bigl(\|\partial_yW\|_{L^2}^2\bigr)^{\frac14}\\
\leq{} & C_m\|u_0\|_{L^2}^{\frac12}\|W\|_{L^2}
\left(\|\partial_yW\|_{L^2}^2
+\frac1{1+t}\int_{\R^2}
V_m\!\left(\frac{y}{\sqrt{1+t}}\right)|W|^2\,dx\,dy\right)^{\frac34}.
\end{align*}
To recover a full power of the quadratic form, we use
\eqref{eq:moving-coercivity} for one half-power of the remaining
weighted norm:
\begin{equation*}
\|W\|_{L^2}^{\frac12}
\leq C_m(1+t)^{\frac14}
\left(\|\partial_yW\|_{L^2}^2
+\frac1{1+t}\int_{\R^2}
V_m\!\left(\frac{y}{\sqrt{1+t}}\right)|W|^2\,dx\,dy\right)^{\frac14}.
\end{equation*}
Splitting $\|W\|_{L^2}$ into two half-powers and applying this
bound to one of them yields
\begin{align}
&\left|\int_{\R^2}
\frac{myu^2}{2(2m-1)(1+t)+y^2}|W|^2\,dx\,dy\right|\notag\\
\leq{} & C_m\|u_0\|_{L^2}^{\frac12}(1+t)^{\frac14}
\|W\|_{L^2}^{\frac12}
\left(\|\partial_yW\|_{L^2}^2
+\frac1{1+t}\int_{\R^2}
V_m\!\left(\frac{y}{\sqrt{1+t}}\right)|W|^2\,dx\,dy\right).
\label{eq:NS-weighted-nonlinear}
\end{align}
The quadratic form can therefore absorb the transport term once
its coefficient is small.

Assume on $[0,T]$ that
\begin{equation}\label{eq:NS-weighted-bootstrap}
\sup_{0\leq t\leq T}
(1+t)^{\frac{3m}{4(m+1)}}\|W(t)\|_{L^2}
\leq2\|w_m\omega_0\|_{L^2}.
\end{equation}
The coefficient in the last line of
\eqref{eq:NS-weighted-nonlinear} is bounded by
\begin{equation*}
C_m\|u_0\|_{L^2}^{\frac{1}{2}}
\|w_m\omega_0\|_{L^2}^{\frac{1}{2}}
(1+t)^{-\frac{m-2}{8(m+1)}}
\leq C_m\varepsilon
(1+t)^{-\frac{m-2}{8(m+1)}}.
\end{equation*}
Choose $\varepsilon_m$ so that this coefficient is at most $1/2$.
Equations \eqref{eq:common-weighted-energy} and
\eqref{eq:NS-weighted-nonlinear} then give
\begin{equation*}\begin{aligned}
&\frac12\frac{d}{dt}\|W\|_{L^2}^2
+\left(1-C_m\varepsilon(1+t)^{-\frac{m-2}{8(m+1)}}\right)
\Biggl[\|\partial_yW\|_{L^2}^2\\
&\qquad+\frac1{1+t}\int_{\R^2}
V_m\!\left(\frac{y}{\sqrt{1+t}}\right)|W|^2\,dx\,dy\Biggr]\leq0.
\end{aligned}\end{equation*}
Since the coefficient of the quadratic form is positive,
\eqref{eq:moving-coercivity} gives
\begin{equation*}
\frac12\frac{d}{dt}\|W\|_{L^2}^2
+\frac{3m}{4(m+1)(1+t)}
\left(1-C_m\varepsilon(1+t)^{-\frac{m-2}{8(m+1)}}\right)
\|W\|_{L^2}^2\leq0.
\end{equation*}
Since $m>2$, the correction has the integrable factor
$(1+t)^{-1-(m-2)/[8(m+1)]}$. Direct integration yields
\begin{equation}\label{eq:NS-weighted-bootstrap-improvement}
\|W(t)\|_{L^2}
\leq\|w_m\omega_0\|_{L^2}
(1+t)^{-\frac{3m}{4(m+1)}}\exp(C_m\varepsilon).
\end{equation}
Decreasing $\varepsilon_m$ makes the exponential at most $3/2$,
which improves \eqref{eq:NS-weighted-bootstrap}. By continuity,
\eqref{eq:NS-weighted-bootstrap-improvement} holds throughout the
interval of existence.

By the definition of $W$, we have proved
\begin{equation}\label{eq:weighted-vorticity-decay}
\left\|w_m\!\left(\frac{y}{\sqrt{1+t}}\right)\omega(t)\right\|_{L^2}
\leq C_m\|w_m\omega_0\|_{L^2}
(1+t)^{-\frac{3m}{4(m+1)}}.
\end{equation}
Since $w_m\geq1$, this also gives
\begin{equation}\label{eq:unweighted-vorticity-decay}
\|\omega(t)\|_{L^2}
\leq C_m\|w_m\omega_0\|_{L^2}
(1+t)^{-\frac{3m}{4(m+1)}}.
\end{equation}
This is \eqref{eq:decay-reduction}, and closes the $H^2$ bootstrap.

\subsection{Velocity decay}

Recovering the velocity from the vorticity introduces a singular
Fourier multiplier at the origin. We control it by combining the
vanishing vertical mean with the first vertical moment. Both follow
from oddness and the weighted estimate just obtained.

First, Cauchy--Schwarz in $y$ gives
\begin{align*}
\|\omega(t)\|_{L_x^2L_y^1}^2
&\leq\left(\int_\R w_m^{-2}\!\left(\frac{y}{\sqrt{1+t}}\right)\,dy\right)
\left\|w_m\!\left(\frac{y}{\sqrt{1+t}}\right)\omega(t)\right\|_{L^2}^2\\
&=(1+t)^{\frac12}\left(\int_\R w_m^{-2}(\eta)\,d\eta\right)
\left\|w_m\!\left(\frac{y}{\sqrt{1+t}}\right)\omega(t)\right\|_{L^2}^2.
\end{align*}
The one-dimensional integral is finite. Hence
\begin{equation*}
\|\omega(t)\|_{L_x^2L_y^1}
\leq C_m(1+t)^{\frac14}
\left\|w_m\!\left(\frac{y}{\sqrt{1+t}}\right)\omega(t)\right\|_{L^2}.
\end{equation*}
The vertical integral is therefore well defined for almost every
$x$, and oddness implies
\[
\int_\R\omega(t,x,y)\,dy=0.
\]
Next, the bound $|\eta|\leq C_mw_m(\eta)$ gives
\begin{equation}\label{eq:first-vertical-moment}
\|y\omega(t)\|_{L^2}
\leq C_m(1+t)^{\frac12}
\left\|w_m\!\left(\frac{y}{\sqrt{1+t}}\right)\omega(t)\right\|_{L^2}.
\end{equation}

Let $(k,\ell)$ be the Fourier variables associated with $(x,y)$.
The preceding two bounds imply, for almost every $k$, that
$\widehat\omega(t,k,\cdot)$ belongs to $H^1(\R)$ and has trace zero
at $\ell=0$. For a smooth approximation with this zero trace,
integration by parts on $(0,\infty)$ gives
\begin{align*}
\int_0^\infty\frac{|\widehat\omega(t,k,\ell)|^2}{\ell^2}\,d\ell
&=2\operatorname{Re}\int_0^\infty
\frac{\overline{\widehat\omega(t,k,\ell)}
\partial_\ell\widehat\omega(t,k,\ell)}{\ell}\,d\ell\\
&\leq2\left(\int_0^\infty
\frac{|\widehat\omega(t,k,\ell)|^2}{\ell^2}\,d\ell\right)^{\frac12}
\left(\int_0^\infty
|\partial_\ell\widehat\omega(t,k,\ell)|^2\,d\ell\right)^{\frac12}.
\end{align*}
Applying the same argument on $(-\infty,0)$ and passing to the
limit yields the one-dimensional Hardy bound
\[
\int_\R\frac{|\widehat\omega(t,k,\ell)|^2}{\ell^2}\,d\ell
\leq4\int_\R|\partial_\ell\widehat\omega(t,k,\ell)|^2\,d\ell.
\]
The Biot--Savart law and Plancherel's theorem now imply
\begin{align}
\|u(t)\|_{L^2}^2
&=\int_{\R^2}\frac{|\widehat\omega(t,k,\ell)|^2}{k^2+\ell^2}
\,dk\,d\ell\notag\\
&\leq\int_\R\left(\int_\R
\frac{|\widehat\omega(t,k,\ell)|^2}{\ell^2}\,d\ell\right)\,dk\notag\\
&\leq4\int_{\R^2}|\partial_\ell\widehat\omega(t,k,\ell)|^2\,dk\,d\ell
=4\|y\omega(t)\|_{L^2}^2.
\label{eq:velocity-from-moment}
\end{align}
Combining \eqref{eq:velocity-from-moment} with
\eqref{eq:first-vertical-moment} and
\eqref{eq:weighted-vorticity-decay}, we obtain
\[
\|u(t)\|_{L^2}
\leq C_m\|w_m\omega_0\|_{L^2}
(1+t)^{\frac12-\frac{3m}{4(m+1)}}
=C_m\|w_m\omega_0\|_{L^2}
(1+t)^{-\frac{m-2}{4(m+1)}}.
\]
This proves \eqref{eq:main-velocity-decay}.

\subsection{Proof of the one-dimensional inequality}
\label{subsec:proof-coercivity}

We now prove the lower bound used in
\eqref{eq:moving-coercivity}. Since the integrand is even for an
odd $f$, the full-line estimate reduces to the positive half-line.
Oddness supplies the zero boundary value at the origin.
The inverse weight has the same rational logarithmic derivative
that enters $V_m$. Multiplication by $\eta$ imposes the zero
boundary value, so we consider the positive comparison function
\begin{equation*}
\phi_m(\eta)=\eta w_m^{-1}(\eta),\qquad \eta>0.
\end{equation*}
Its derivatives are
\begin{align*}
\phi_m'
&=w_m^{-1}\left(1-\frac{m\eta^2}{2(2m-1)+\eta^2}\right),
\\*
\phi_m''
&=\phi_m\left(-\frac{3m}{2(2m-1)+\eta^2}
+\frac{m(m+2)\eta^2}{[2(2m-1)+\eta^2]^2}\right).
\end{align*}
Using \eqref{eq:potential}, we find
\begin{equation*}
-\frac{\phi_m''}{\phi_m}+V_m
=\frac{6m(2m-1)+\frac{m}{2}\eta^4}{[2(2m-1)+\eta^2]^2}.
\end{equation*}
The derivative of the expression on the right-hand side is
\[
\frac{4m(2m-1)\eta(\eta^2-6)}{[2(2m-1)+\eta^2]^3}.
\]
Its minimum on $(0,\infty)$ is attained at $\eta=\sqrt6$ and equals
$3m/[4(m+1)]$. Equivalently, completing the square gives
\begin{equation}\label{eq:comparison-pointwise-bound}
-\frac{\phi_m''}{\phi_m}+V_m-\frac{3m}{4(m+1)}
=\frac{m(2m-1)(\eta^2-6)^2}
{4(m+1)[2(2m-1)+\eta^2]^2}\geq0.
\end{equation}

To turn this pointwise bound into an integral inequality, we
factor the comparison function out of an arbitrary test function.
For $f\in C_c^\infty(0,\infty)$, write $f=\phi_m h$. The factor
$h$ measures the variation relative to $\phi_m$; after integration
by parts, its derivatives will contribute a nonnegative term.
The product rule gives
\[
f''=\phi_m''h+2\phi_m'h'+\phi_m h''.
\]
Multiplying $-f''+V_mf$ by $f$ and integrating, we obtain
\begin{align*}
\int_0^\infty(|f'|^2+V_m f^2)\,d\eta
&=\int_0^\infty(-f''+V_m f)f\,d\eta
\\*
&=\int_0^\infty(-\phi_m''+V_m\phi_m)\phi_m h^2\,d\eta
-2\int_0^\infty\phi_m\phi_m'h'h\,d\eta
\\*
&\quad-\int_0^\infty\phi_m^2h''h\,d\eta.
\end{align*}
The last integral satisfies
\[
-\int_0^\infty\phi_m^2h''h\,d\eta
=\int_0^\infty\phi_m^2|h'|^2\,d\eta
+2\int_0^\infty\phi_m\phi_m'h'h\,d\eta.
\]
The mixed terms cancel. Hence
\begin{equation}\label{eq:comparison-identity}
\int_0^\infty(|f'|^2+V_m f^2)\,d\eta
=\int_0^\infty\phi_m^2
\left|\left(\frac f{\phi_m}\right)'\right|^2\,d\eta
+\int_0^\infty\left(-\frac{\phi_m''}{\phi_m}+V_m\right)f^2\,d\eta.
\end{equation}
The first term is nonnegative. Applying
\eqref{eq:comparison-pointwise-bound} proves the half-line
inequality. Density in $H_0^1(0,\infty)$ and reflection prove
Proposition~\ref{prop:coercive-inequality}.

The Gaussian-weight assertion in Remark~\ref{rem:gaussian} follows
from the same calculation. For the weight $e^{\eta^2/8}$, the
potential is $\eta^2/16$, and the positive comparison function
$\eta e^{-\eta^2/8}$ satisfies
\[
-\frac{(\eta e^{-\frac{\eta^2}{8}})''}{\eta e^{-\frac{\eta^2}{8}}}
+\frac{\eta^2}{16}=\frac34.
\]
Thus \eqref{eq:comparison-identity} gives the odd-function lower
bound with constant $3/4$. The Navier--Stokes transport estimate
\eqref{eq:NS-weighted-holder} remains valid, since the square of
the logarithmic derivative of the Gaussian weight equals its
potential. The bootstrap above therefore gives
\[
\left\|e^{\frac{y^2}{8(1+t)}}\omega(t)\right\|_{L^2}
\leq C\|e^{\frac{y^2}{8}}\omega_0\|_{L^2}(1+t)^{-\frac34}.
\]
The moment estimate \eqref{eq:velocity-from-moment} then yields
velocity decay with exponent $1/4$.

\subsection{Higher-order stability}

We prove the higher-regularity assertion of Theorem~\ref{thm:main}
under \eqref{eq:H3-small-data}. For $T>0$, define
\begin{equation}\label{eq:E3-definition}
	E_3(T)
	=
	\sup_{0\leq t\leq T}\|u(t)\|_{H^3}^2
	+
	\int_0^T\|\partial_yu(t)\|_{H^3}^2\,dt.
\end{equation}
Fix $C_0\geq1$, depending only on $m$, to dominate the combined constants
multiplying the initial energy and the established $H^2$ bound after
the norm comparisons below. Constants depending on this fixed choice
will be denoted by $C_m$. Assume that
\begin{equation}\label{eq:E3-bootstrap}
	E_3(T)\leq4C_0\varepsilon^2.
\end{equation}

The $H^2$ estimate established above gives
\begin{equation}\label{eq:H2-bound-for-H3}
	\sup_{0\leq t\leq T}\|u(t)\|_{H^2}^2
	+
	\int_0^T\|\partial_yu(t)\|_{H^2}^2\,dt
	\leq
	C_m\varepsilon^2.
\end{equation}
It therefore remains to control only the homogeneous highest-order
part of \eqref{eq:E3-definition}.

By the Fourier characterization of homogeneous Sobolev norms and the
Biot--Savart law,
\begin{align}
	\|u\|_{\dot H^3}^2
	&\approx
	\|\partial_x^3u\|_{L^2}^2
	+
	\|\partial_y^3u\|_{L^2}^2
	\approx
	\|\partial_x^2\omega\|_{L^2}^2
	+
	\|\partial_y^2\omega\|_{L^2}^2,
	\label{eq:H3-pure-derivative-equivalence}\\*
	\|\partial_yu\|_{\dot H^3}^2
	&\approx
	\|\partial_{xxxy}u\|_{L^2}^2
	+
	\|\partial_y^4u\|_{L^2}^2
	\approx
	\|\partial_{xxy}\omega\|_{L^2}^2
	+
	\|\partial_y^3\omega\|_{L^2}^2.
	\label{eq:H3-dissipation-equivalence}
\end{align}
Consequently,
\begin{align}
	E_3(T)
	&\approx
	E_2(T)
	+
	\sup_{0\leq t\leq T}
	\left(
	\|\partial_x^2\omega(t)\|_{L^2}^2
	+
	\|\partial_y^2\omega(t)\|_{L^2}^2
	\right)
	\notag\\*
	&\quad+
	\int_0^T
	\left(
	\|\partial_{xxy}\omega(t)\|_{L^2}^2
	+
	\|\partial_y^3\omega(t)\|_{L^2}^2
	\right)\,dt.
	\label{eq:E3-vorticity-reduction}
\end{align}
Thus it suffices to estimate $\partial_x^2\omega$ and
$\partial_y^2\omega$.

We begin with the second-order horizontal derivative.
Applying $\partial_x^2$ to \eqref{eq:vorticity-equation} and taking the
$L^2$ inner product with $\partial_x^2\omega$, we obtain
\begin{equation}\begin{aligned}
	\frac12\frac{d}{dt}\|\partial_x^2\omega\|_{L^2}^2
	+\|\partial_{xxy}\omega\|_{L^2}^2
	&=
	-\int_{\R^2}
	\partial_x^2u^1\,
	\partial_x\omega\,
	\partial_x^2\omega\,dx\,dy
	-\int_{\R^2}
	\partial_x^2u^2\,
	\partial_y\omega\,
	\partial_x^2\omega\,dx\,dy
	\\
	&\quad
	-2\int_{\R^2}
	\partial_xu^1
	(\partial_x^2\omega)^2\,dx\,dy
	-2\int_{\R^2}
	\partial_xu^2\,
	\partial_{xy}\omega\,
	\partial_x^2\omega\,dx\,dy
	\\
	&=:B_1+B_2+B_3+B_4.
\end{aligned}\label{eq:horizontal-second-energy}\end{equation}

For $B_1$, incompressibility and integration by parts in $y$ give
\begin{align}
	B_1
	&=
	-\int_{\R^2}
	\partial_xu^2\,
	\partial_{xy}\omega\,
	\partial_x^2\omega\,dx\,dy
	-
	\int_{\R^2}
	\partial_xu^2\,
	\partial_x\omega\,
	\partial_{xxy}\omega\,dx\,dy.
	\label{eq:B1-splitting}
\end{align}
For the second integral in \eqref{eq:B1-splitting}, the mixed
derivative is estimated by
\begin{equation}\label{eq:H3-interpolation}
\|\partial_{xy}\omega\|_{L^2}
\leq\|\partial_x\omega\|_{L^2}^{\frac12}
\|\partial_{xyy}\omega\|_{L^2}^{\frac12}
\leq C\|\omega\|_{L^2}^{\frac14}
\|\omega\|_{\dot H^2}^{\frac14}
\|\partial_y\omega\|_{\dot H^2}^{\frac12}.
\end{equation}
Anisotropic H\"older's inequality, one-dimensional interpolation,
and Young's inequality now give
\begin{align}
	|B_1|
	&\leq
	\|\partial_xu^2\|_{L^2}
	\|\partial_{xy}\omega\|_{L_y^2L_x^\infty}
	\|\partial_x^2\omega\|_{L_x^2L_y^\infty}
	\notag\\
	&\quad+
	\|\partial_xu^2\|_{L_y^2L_x^\infty}
	\|\partial_x\omega\|_{L_x^2L_y^\infty}
	\|\partial_{xxy}\omega\|_{L^2}
	\notag\\
	&\leq C\|\omega\|_{L^2}
	\|\partial_{xy}\omega\|_{L^2}^{\frac{1}{2}}
	\|\partial_x^2\omega\|_{L^2}^{\frac{1}{2}}
	\|\partial_{xxy}\omega\|_{L^2}\notag\\
	&\quad+C\|\omega\|_{L^2}^{\frac{1}{2}}
	\|\partial_x\omega\|_{L^2}
	\|\partial_{xy}\omega\|_{L^2}^{\frac{1}{2}}
	\|\partial_{xxy}\omega\|_{L^2}\notag\\
	&\leq C\|\omega\|_{L^2}\|\omega\|_{\dot H^2}
	\|\partial_y\omega\|_{\dot H^2}
	+C\|\omega\|_{L^2}^{\frac{9}{8}}\|\omega\|_{\dot H^2}^{\frac{5}{8}}
	\|\partial_y\omega\|_{\dot H^2}^{\frac{5}{4}}\notag\\
	&\leq
	2\delta\|\partial_y\omega\|_{\dot H^2}^2
	+
	C_\delta\|\omega\|_{L^2}^2
	\|\omega\|_{\dot H^2}^2+
	C_\delta\|\omega\|_{L^2}^3
	\|\omega\|_{\dot H^2}^{\frac{5}{3}}.
	\label{eq:B1-estimate}
\end{align}

For $B_2$, integration by parts in $y$ yields
\begin{align}
	B_2
	&=
	\int_{\R^2}
	\partial_{xxy}u^2\,
	\omega\,
	\partial_x^2\omega\,dx\,dy
	+
	\int_{\R^2}
	\partial_x^2u^2\,
	\omega\,
	\partial_{xxy}\omega\,dx\,dy.
	\label{eq:B2-splitting}
\end{align}
For the first integral, we use the Biot--Savart bound and
\eqref{eq:H3-interpolation}. For the second, interpolation in $y$
gives the additional estimate
\[
\|\partial_y\omega\|_{L^2}
\leq\|\omega\|_{L^2}^{\frac23}
\|\partial_y^3\omega\|_{L^2}^{\frac13}
\leq\|\omega\|_{L^2}^{\frac23}
\|\partial_y\omega\|_{\dot H^2}^{\frac13}.
\]
Consequently,
\begin{align}
	|B_2|
	&\leq
	\|\partial_{xxy}u^2\|_{L^2}
	\|\omega\|_{L_y^2L_x^\infty}
	\|\partial_x^2\omega\|_{L_x^2L_y^\infty}
	\notag\\
	&\quad+
	\|\partial_x^2u^2\|_{L_y^2L_x^\infty}
	\|\omega\|_{L_x^2L_y^\infty}
	\|\partial_{xxy}\omega\|_{L^2}
	\notag\\
	&\leq C\|\partial_{xy}\omega\|_{L^2}
	\|\omega\|_{L^2}^{\frac{1}{2}}
	\|\partial_x\omega\|_{L^2}^{\frac{1}{2}}
	\|\partial_x^2\omega\|_{L^2}^{\frac{1}{2}}
	\|\partial_{xxy}\omega\|_{L^2}^{\frac{1}{2}}\notag\\
	&\quad+C\|\partial_x\omega\|_{L^2}^{\frac{1}{2}}
	\|\partial_x^2\omega\|_{L^2}^{\frac{1}{2}}
	\|\omega\|_{L^2}^{\frac{1}{2}}
	\|\partial_y\omega\|_{L^2}^{\frac{1}{2}}
	\|\partial_{xxy}\omega\|_{L^2}\notag\\
	&\leq C\|\omega\|_{L^2}\|\omega\|_{\dot H^2}
	\|\partial_y\omega\|_{\dot H^2}
	+C\|\omega\|_{L^2}^{\frac{13}{12}}\|\omega\|_{\dot H^2}^{\frac{3}{4}}
	\|\partial_y\omega\|_{\dot H^2}^{\frac{7}{6}}\notag\\
	&\leq
	2\delta\|\partial_y\omega\|_{\dot H^2}^2
	+
	C_\delta\|\omega\|_{L^2}^2
	\|\omega\|_{\dot H^2}^2+
	C_\delta\|\omega\|_{L^2}^{\frac{13}{5}}
	\|\omega\|_{\dot H^2}^{\frac{9}{5}}.
	\label{eq:B2-estimate}
\end{align}

For $B_3$, incompressibility and integration by parts in $y$ give
\[
B_3
=
-4\int_{\R^2}
u^2\,
\partial_{xxy}\omega\,
\partial_x^2\omega\,dx\,dy.
\]
Therefore,
\begin{align}
	|B_3|
	&\leq
	4\|u^2\|_{L_y^2L_x^\infty}
	\|\partial_x^2\omega\|_{L_x^2L_y^\infty}
	\|\partial_{xxy}\omega\|_{L^2}
	\notag\\*
	&\leq
	C\|u^2\|_{L^2}^{\frac{1}{2}}
	\|\partial_xu^2\|_{L^2}^{\frac{1}{2}}
	\|\partial_x^2\omega\|_{L^2}^{\frac{1}{2}}
	\|\partial_{xxy}\omega\|_{L^2}^{\frac{3}{2}}
	\notag\\*
	&\leq
	\delta\|\partial_y\omega\|_{\dot H^2}^2
	+
	C_\delta\|u_0\|_{L^2}^2
	\|\omega\|_{L^2}^2
	\|\omega\|_{\dot H^2}^2.
	\label{eq:B3-estimate}
\end{align}

Finally, $B_4$ is estimated in the same way as the first term in
\eqref{eq:B1-splitting}:
\begin{align}
	|B_4|
	&\leq
	C\|\partial_xu^2\|_{L^2}
	\|\partial_{xy}\omega\|_{L_y^2L_x^\infty}
	\|\partial_x^2\omega\|_{L_x^2L_y^\infty}
	\notag\\*
	&\leq
	\delta\|\partial_y\omega\|_{\dot H^2}^2
	+
	C_\delta\|\omega\|_{L^2}^2
	\|\omega\|_{\dot H^2}^2.
	\label{eq:B4-estimate}
\end{align}

Collecting \eqref{eq:B1-estimate}--\eqref{eq:B4-estimate}, we find
\begin{align}
	\sum_{j=1}^4|B_j|
	&\leq
	C\delta\|\partial_y\omega\|_{\dot H^2}^2
	+
	C_\delta\bigl(1+\|u_0\|_{L^2}^2\bigr)
	\|\omega\|_{L^2}^2
	\|\omega\|_{\dot H^2}^2
	\notag\\*
	&\quad+
	C_\delta\|\omega\|_{L^2}^3
	\|\omega\|_{\dot H^2}^{\frac{5}{3}}
	+
	C_\delta\|\omega\|_{L^2}^{\frac{13}{5}}
	\|\omega\|_{\dot H^2}^{\frac{9}{5}}.
	\label{eq:B-sum-estimate}
\end{align}

We next estimate the second-order vertical derivative. Applying $\partial_y^2$ to \eqref{eq:vorticity-equation} and taking the
$L^2$ inner product with $\partial_y^2\omega$, we obtain
\begin{align}
	\frac12\frac{d}{dt}\|\partial_y^2\omega\|_{L^2}^2
	+\|\partial_y^3\omega\|_{L^2}^2
	&=
	-\int_{\R^2}
	\partial_y^2u^1\,
	\partial_x\omega\,
	\partial_y^2\omega\,dx\,dy
	-\int_{\R^2}
	\partial_y^2u^2\,
	\partial_y\omega\,
	\partial_y^2\omega\,dx\,dy
	\notag\\*
	&\quad
	-2\int_{\R^2}
	\partial_yu^1\,
	\partial_{xy}\omega\,
	\partial_y^2\omega\,dx\,dy
	-2\int_{\R^2}
	\partial_yu^2
	(\partial_y^2\omega)^2\,dx\,dy
	\notag\\*
	&=:D_1+D_2+D_3+D_4.
	\label{eq:vertical-second-energy}
\end{align}

For $D_1$, anisotropic H\"older's inequality gives
\begin{align*}
	|D_1|
	&\leq
	\|\partial_x\omega\|_{L^2}
	\|\partial_y^2u^1\|_{L_y^2L_x^\infty}
	\|\partial_y^2\omega\|_{L_x^2L_y^\infty}
	\\*
	&\leq
	C\|\partial_x\omega\|_{L^2}
	\|\partial_y^2u^1\|_{L^2}^{\frac{1}{2}}
	\|\partial_{xyy}u^1\|_{L^2}^{\frac{1}{2}}
	\|\partial_y^2\omega\|_{L^2}^{\frac{1}{2}}
	\|\partial_y^3\omega\|_{L^2}^{\frac{1}{2}}
	\\*
	&\leq
	C\|\omega\|_{H^2}
	\|\partial_y\omega\|_{H^2}^2.
\end{align*}

Similarly,
\begin{align*}
	|D_2|
	&\leq
	\|\partial_y\omega\|_{L^2}
	\|\partial_y^2u^2\|_{L_y^2L_x^\infty}
	\|\partial_y^2\omega\|_{L_x^2L_y^\infty}
	\\*
	&\leq
	C\|\partial_y\omega\|_{L^2}
	\|\partial_y^2u^2\|_{L^2}^{\frac{1}{2}}
	\|\partial_{xyy}u^2\|_{L^2}^{\frac{1}{2}}
	\|\partial_y^2\omega\|_{L^2}^{\frac{1}{2}}
	\|\partial_y^3\omega\|_{L^2}^{\frac{1}{2}}
	\\*
	&\leq
	C\|\omega\|_{H^2}
	\|\partial_y\omega\|_{H^2}^2.
\end{align*}

For $D_3$, we use
\[
\partial_{xy}u^1=-\partial_y^2u^2,
\qquad
\|\partial_{xy}u^1\|_{L^2}
\leq
\|\partial_y\omega\|_{L^2},
\]
and obtain
\begin{align*}
	|D_3|
	&\leq
	C\|\partial_{xy}\omega\|_{L^2}
	\|\partial_yu^1\|_{L_y^2L_x^\infty}
	\|\partial_y^2\omega\|_{L_x^2L_y^\infty}
	\\*
	&\leq
	C\|\partial_{xy}\omega\|_{L^2}
	\|\partial_yu^1\|_{L^2}^{\frac{1}{2}}
	\|\partial_{xy}u^1\|_{L^2}^{\frac{1}{2}}
	\|\partial_y^2\omega\|_{L^2}^{\frac{1}{2}}
	\|\partial_y^3\omega\|_{L^2}^{\frac{1}{2}}
	\\*
	&\leq
	C\|\omega\|_{H^2}
	\|\partial_y\omega\|_{H^2}^2.
\end{align*}
Here we also used
\[
\|\partial_{xy}\omega\|_{L^2}
\leq C\|\omega\|_{H^2}^{\frac12}\|\partial_y\omega\|_{H^2}^{\frac12}.
\]

Finally,
\begin{align*}
	|D_4|
	&\leq
	C\|\partial_yu^2\|_{L^2}
	\|\partial_y^2\omega\|_{L_y^2L_x^\infty}
	\|\partial_y^2\omega\|_{L_x^2L_y^\infty}
	\\*
	&\leq
	C\|\partial_yu^2\|_{L^2}
	\|\partial_y^2\omega\|_{L^2}
	\|\partial_{xyy}\omega\|_{L^2}^{\frac{1}{2}}
	\|\partial_y^3\omega\|_{L^2}^{\frac{1}{2}}
	\\*
	&\leq
	C\|\omega\|_{H^2}
	\|\partial_y\omega\|_{H^2}^2.
\end{align*}
It follows that
\begin{equation}\label{eq:D-sum-estimate}
	\sum_{j=1}^4|D_j|
	\leq
	C\|\omega\|_{H^2}
	\|\partial_y\omega\|_{H^2}^2.
\end{equation}

We now combine the horizontal and vertical estimates. By
\eqref{eq:H2-equivalence}, \eqref{eq:H3-pure-derivative-equivalence},
and \eqref{eq:E3-bootstrap},
\begin{equation}\label{eq:H2-vorticity-bootstrap}
	\sup_{0\leq t\leq T}\|\omega(t)\|_{H^2}
	\leq
	C_m\varepsilon.
\end{equation}
Moreover, \eqref{eq:H2-bound-for-H3} gives
\begin{equation}\label{eq:lower-order-dissipation-H3}
	\int_0^T
	\|\partial_y\omega(t)\|_{H^1}^2\,dt
	\leq
	C_m\varepsilon^2.
\end{equation}
Therefore, \eqref{eq:D-sum-estimate} yields
\begin{align}
	\sum_{j=1}^4|D_j|
	&\leq
	C_m\varepsilon
	\|\partial_y\omega\|_{\dot H^2}^2
	+
	C_m\varepsilon
	\|\partial_y\omega\|_{H^1}^2.
	\label{eq:D-absorbed-estimate}
\end{align}

We first choose $\delta>0$ sufficiently small in
\eqref{eq:B-sum-estimate} and then choose $\varepsilon_m>0$
sufficiently small in \eqref{eq:D-absorbed-estimate}. Adding
\eqref{eq:horizontal-second-energy} and
\eqref{eq:vertical-second-energy}, we obtain
\begin{align}
	\frac{d}{dt}\bigl(\|\partial_x^2\omega\|_{L^2}^2+\|\partial_y^2\omega\|_{L^2}^2\bigr)
	+
	\frac12\|\partial_y\omega\|_{\dot H^2}^2
	&\leq
	C\bigl(1+\|u_0\|_{L^2}^2\bigr)
	\|\omega\|_{L^2}^2
	\|\omega\|_{\dot H^2}^2
	\notag\\*
	&\quad+
	C\|\omega\|_{L^2}^3
	\|\omega\|_{\dot H^2}^{\frac{5}{3}}
	+
	C\|\omega\|_{L^2}^{\frac{13}{5}}
	\|\omega\|_{\dot H^2}^{\frac{9}{5}}
	\notag\\*
	&\quad+
	C\|\omega\|_{H^2}
	\|\partial_y\omega\|_{H^1}^2.
	\label{eq:H3-homogeneous-differential}
\end{align}

Since $m>2$, the decay estimate
\eqref{eq:unweighted-vorticity-decay} and the bootstrap assumption give
\begin{align}
	\int_0^T
	\Bigl[
	\bigl(1+\|u_0\|_{L^2}^2\bigr)
	\|\omega\|_{L^2}^2
	\|\omega\|_{\dot H^2}^2
	+
	\|\omega\|_{L^2}^3
	\|\omega\|_{\dot H^2}^{\frac{5}{3}}
	+
	\|\omega\|_{L^2}^{\frac{13}{5}}
	\|\omega\|_{\dot H^2}^{\frac{9}{5}}
	\Bigr]\,dt
	&\leq
	C_m\varepsilon^4.
	\label{eq:B-time-integral}
\end{align}
On the other hand,
\begin{equation}\label{eq:D-time-integral}
	\int_0^T
	\|\omega\|_{H^2}
	\|\partial_y\omega\|_{H^1}^2\,dt
	\leq
	C_m\varepsilon^3
\end{equation}
by \eqref{eq:H2-vorticity-bootstrap} and
\eqref{eq:lower-order-dissipation-H3}.

Integrating \eqref{eq:H3-homogeneous-differential}, using
\eqref{eq:H3-pure-derivative-equivalence}--\eqref{eq:H3-dissipation-equivalence}
and
\eqref{eq:B-time-integral}--\eqref{eq:D-time-integral}, we find
\begin{align}
	\sup_{0\leq t\leq T}\|\omega(t)\|_{\dot H^2}^2
	+
	\int_0^T
	\|\partial_y\omega(t)\|_{\dot H^2}^2\,dt
&\leq
	C\|u_0\|_{H^3}^2
	+
	C_m\varepsilon^3.
	\label{eq:H3-homogeneous-estimate}
\end{align}
Combining \eqref{eq:H3-homogeneous-estimate} with
\eqref{eq:H2-bound-for-H3} and
\eqref{eq:E3-vorticity-reduction}, we conclude that
\begin{equation*}
	E_3(T)
	\leq
	C_0\varepsilon^2
	+
	C_m\varepsilon^3.
\end{equation*}
With $C_0$ already fixed, choose $\varepsilon_m$ so that
$C_m\varepsilon_m\leq C_0$. This improves \eqref{eq:E3-bootstrap} to
\[
E_3(T)\leq2C_0\varepsilon^2.
\]
The standard continuity argument closes the bootstrap and proves
\eqref{eq:main-H3-estimate}.

\begin{proof}[Proof of Corollary~\ref{cor:half-plane}]
Extend $u_0^1$ evenly and $u_0^2$ oddly across $y=0$. The boundary
conditions give a divergence free $H^2(\R^2)$ extension with
\eqref{eq:initial-symmetry}. For $H^3$ data, the additional trace
required for the odd extension follows from
\[
\partial_y^2u_0^2=-\partial_x\partial_yu_0^1=0
\qquad\text{on }y=0.
\]
Apply Theorem~\ref{thm:main} to the extension and
restrict to $y>0$, adjusting $\varepsilon_m$ by the fixed extension
factor. The preserved symmetry gives \eqref{eq:half-plane-boundary},
and all the stated estimates follow by restriction.
\end{proof}

\section{Stability and decay for IPM}\label{sec:IPM}

In this section, we prove Theorem~\ref{thm:IPM} by coupling the
density energy estimate to its weighted decay. Darcy's law provides directional identities for the
highest-order horizontal interactions and a decomposition of the
vertical velocity for the weighted estimate.
The calculations are performed for smooth solutions and extend to
$H^2$ data by symmetry-preserving regularization and the uniform
bounds below. Uniqueness follows from an $L^2$ difference estimate
using the directional multipliers in \eqref{eq:IPM-multipliers}.

\subsection{The \texorpdfstring{$H^2$}{H2} energy estimate}

Taking the Fourier transform of Darcy's law gives
\begin{equation}\label{eq:IPM-multipliers}
\widehat{u^1}(k,\ell)
=\frac{k\ell}{k^2+\ell^2}\widehat\rho(k,\ell),
\qquad
\widehat{u^2}(k,\ell)
=-\frac{k^2}{k^2+\ell^2}\widehat\rho(k,\ell).
\end{equation}
In particular,
\begin{equation}\label{eq:IPM-L2-multiplier}
\|\partial_x^a\partial_y^b u\|_{L^2}
\leq\|\partial_x^a\partial_y^b\rho\|_{L^2},
\qquad a,b\in\{0,1,2,\ldots\}.
\end{equation}
We shall also use the identities
\begin{equation}\label{eq:IPM-Darcy-identities}
\partial_xu^1=-\partial_yu^2,
\qquad
\partial_xu^2=\partial_yu^1-\partial_x\rho,
\qquad
\partial_x^2u^2=\partial_{xy}u^1-\partial_x^2\rho.
\end{equation}
The first multiplier in \eqref{eq:IPM-multipliers} yields the more
specific bounds
\begin{align}
\|\partial_{xy}u^1\|_{L^2}
&\leq\|\partial_y^2\rho\|_{L^2},
&\|\partial_x^3u^1\|_{L^2}
&\leq\|\partial_{xxy}\rho\|_{L^2},\notag\\
\|\partial_{xxy}u^1\|_{L^2}
&\leq\|\partial_{xxy}\rho\|_{L^2},
&\|\partial_{xyy}u^1\|_{L^2}
&\leq\|\partial_y^3\rho\|_{L^2}.
\label{eq:IPM-directional-multipliers}
\end{align}
These estimates express the vertical factor in the multiplier of
$u^1$, which will be used below to control horizontal interactions.

Taking the $L^2$ inner product of the density equation with $\rho$
gives
\begin{equation}\label{eq:IPM-density-L2}
\|\rho(t)\|_{L^2}^2
+2\int_0^t\|\partial_y\rho(r)\|_{L^2}^2\,dr
=\|\rho_0\|_{L^2}^2.
\end{equation}
Since
\begin{align}
\|\rho\|_{H^2}^2
&\approx\|\rho\|_{L^2}^2
+\|\partial_x^2\rho\|_{L^2}^2+\|\partial_y^2\rho\|_{L^2}^2,
\notag\\*
\|\partial_y\rho\|_{H^2}^2
&\approx\|\partial_y\rho\|_{L^2}^2
+\|\partial_{xxy}\rho\|_{L^2}^2+\|\partial_y^3\rho\|_{L^2}^2,
\label{eq:IPM-norm-equivalences}
\end{align}
it remains to estimate the two pure second derivatives.

Applying $\partial_x^2$ to the density equation and taking the
$L^2$ inner product with $\partial_x^2\rho$, we obtain
\begin{align*}
\frac12\frac{d}{dt}\|\partial_x^2\rho\|_{L^2}^2
+\|\partial_{xxy}\rho\|_{L^2}^2
&=-\int_{\R^2}\partial_x^2u^1\partial_x\rho\partial_x^2\rho\,dx\,dy
-\int_{\R^2}\partial_x^2u^2\partial_y\rho\partial_x^2\rho\,dx\,dy
\\*
&\quad-2\int_{\R^2}\partial_xu^1(\partial_x^2\rho)^2\,dx\,dy
-2\int_{\R^2}\partial_xu^2\partial_{xy}\rho\partial_x^2\rho\,dx\,dy
\\*
&=:G_1+G_2+G_3+G_4.
\end{align*}

For $G_1$, integration by parts in $x$ gives
\[
G_1=\frac12\int_{\R^2}\partial_x^3u^1(\partial_x\rho)^2\,dx\,dy.
\]
The directional multiplier estimate in
\eqref{eq:IPM-directional-multipliers} controls the velocity factor.
To estimate the remaining density factors, we use
\begin{align}
\|\partial_x\rho\|_{L^2}
&\leq\|\rho\|_{L^2}^{\frac12}\|\rho\|_{\dot H^2}^{\frac12},
\notag\\
\|\partial_{xy}\rho\|_{L^2}
&\leq C\|\rho\|_{L^2}^{\frac13}
\|\partial_y\rho\|_{\dot H^2}^{\frac23}.
\label{eq:IPM-interpolation}
\end{align}
The second inequality follows by H\"older's inequality in Fourier
space from
\[
|k\ell|^{\frac32}\leq C|\ell|(k^2+\ell^2).
\]
Anisotropic H\"older's inequality and one-dimensional interpolation
therefore yield
\begin{align}
|G_1|
&\leq C\|\partial_x^3u^1\|_{L^2}
\|\partial_x\rho\|_{L_y^2L_x^\infty}
\|\partial_x\rho\|_{L_x^2L_y^\infty}
\notag\\*
&\leq C\|\partial_{xxy}\rho\|_{L^2}
\|\partial_x\rho\|_{L^2}
\|\partial_x^2\rho\|_{L^2}^{\frac{1}{2}}
\|\partial_{xy}\rho\|_{L^2}^{\frac{1}{2}}
\notag\\*
&\leq C\|\rho\|_{L^2}^{\frac{2}{3}}\|\rho\|_{\dot H^2}
\|\partial_y\rho\|_{\dot H^2}^{\frac{4}{3}}.
\label{eq:IPM-B1}
\end{align}

For $G_2$, the last identity in
\eqref{eq:IPM-Darcy-identities} and integration by parts in $y$
give
\begin{align*}
G_2
&=-\int_{\R^2}\partial_{xy}u^1\partial_y\rho\partial_x^2\rho\,dx\,dy
+\int_{\R^2}\partial_y\rho(\partial_x^2\rho)^2\,dx\,dy\\*
&=-\int_{\R^2}\partial_{xy}u^1\partial_y\rho\partial_x^2\rho\,dx\,dy
-2\int_{\R^2}\rho\partial_x^2\rho\partial_{xxy}\rho\,dx\,dy.
\end{align*}
For this estimate, interpolation in the vertical variable gives
\begin{equation}\begin{aligned}
\|\partial_y\rho\|_{L^2}
&\leq\|\rho\|_{L^2}^{\frac23}
\|\partial_y^3\rho\|_{L^2}^{\frac13},\\
\|\partial_y^2\rho\|_{L^2}
&\leq\|\rho\|_{L^2}^{\frac13}
\|\partial_y^3\rho\|_{L^2}^{\frac23}.
\end{aligned}\label{eq:IPM-vertical-interpolation}\end{equation}
Together with \eqref{eq:IPM-interpolation}, this yields
\begin{align}
|G_2|
&\leq
\|\partial_{xy}u^1\|_{L_y^2L_x^\infty}
\|\partial_y\rho\|_{L_x^2L_y^\infty}
\|\partial_x^2\rho\|_{L^2}
\notag\\
&\quad+C\|\rho\|_{L_y^2L_x^\infty}
\|\partial_x^2\rho\|_{L_x^2L_y^\infty}
\|\partial_{xxy}\rho\|_{L^2}
\notag\\
&\leq C\|\partial_y^2\rho\|_{L^2}^{\frac{1}{2}}
\|\partial_{xxy}\rho\|_{L^2}^{\frac{1}{2}}
\|\partial_y\rho\|_{L^2}^{\frac{1}{2}}
\|\partial_y^2\rho\|_{L^2}^{\frac{1}{2}}
\|\partial_x^2\rho\|_{L^2}
\notag\\
&\quad+C\|\rho\|_{L^2}^{\frac{1}{2}}\|\partial_x\rho\|_{L^2}^{\frac{1}{2}}
\|\partial_x^2\rho\|_{L^2}^{\frac{1}{2}}
\|\partial_{xxy}\rho\|_{L^2}^{\frac{3}{2}}
\notag\\
&\leq C\|\rho\|_{L^2}^{\frac{2}{3}}\|\rho\|_{\dot H^2}
\|\partial_y\rho\|_{\dot H^2}^{\frac{4}{3}}
+C\|\rho\|_{L^2}^{\frac{3}{4}}\|\rho\|_{\dot H^2}^{\frac{3}{4}}
\|\partial_y\rho\|_{\dot H^2}^{\frac{3}{2}}.
\label{eq:IPM-B2}
\end{align}

For $G_3$, incompressibility and integration by parts in $y$ give
\[
G_3=-4\int_{\R^2}u^2\partial_x^2\rho\partial_{xxy}\rho\,dx\,dy.
\]
The multiplier bounds imply
\[
\|u^2\|_{L^2}\leq\|\rho\|_{L^2},
\qquad
\|\partial_xu^2\|_{L^2}\leq\|\partial_x\rho\|_{L^2}.
\]
Applying the estimate for the second integral in $G_2$, we obtain
\begin{align}
|G_3|
&\leq C\|u^2\|_{L_y^2L_x^\infty}
\|\partial_x^2\rho\|_{L_x^2L_y^\infty}
\|\partial_{xxy}\rho\|_{L^2}
\notag\\*
&\leq C\|\rho\|_{L^2}^{\frac{3}{4}}\|\rho\|_{\dot H^2}^{\frac{3}{4}}
\|\partial_y\rho\|_{\dot H^2}^{\frac{3}{2}}.
\label{eq:IPM-B3}
\end{align}

Finally, the second identity in \eqref{eq:IPM-Darcy-identities}
gives
\begin{align*}
G_4
&=-2\int_{\R^2}\partial_yu^1\partial_{xy}\rho\partial_x^2\rho\,dx\,dy
+2\int_{\R^2}\partial_x\rho\partial_{xy}\rho\partial_x^2\rho\,dx\,dy\\*
&=-2\int_{\R^2}\partial_yu^1\partial_{xy}\rho\partial_x^2\rho\,dx\,dy
-\int_{\R^2}(\partial_x\rho)^2\partial_{xxy}\rho\,dx\,dy.
\end{align*}
The last integral is estimated as in \eqref{eq:IPM-B1}. For the
first integral, we place the two density factors in the opposite
mixed norms and obtain
\begin{align}
|G_4|
&\leq C\|\partial_yu^1\|_{L^2}
\|\partial_{xy}\rho\|_{L_y^2L_x^\infty}
\|\partial_x^2\rho\|_{L_x^2L_y^\infty}
+C\|\rho\|_{L^2}^{\frac{2}{3}}\|\rho\|_{\dot H^2}
\|\partial_y\rho\|_{\dot H^2}^{\frac{4}{3}}
\notag\\
&\leq C\|\partial_y\rho\|_{L^2}
\|\partial_{xy}\rho\|_{L^2}^{\frac{1}{2}}
\|\partial_{xxy}\rho\|_{L^2}^{\frac{1}{2}}
\|\partial_x^2\rho\|_{L^2}^{\frac{1}{2}}
\|\partial_{xxy}\rho\|_{L^2}^{\frac{1}{2}}
\notag\\
&\quad+C\|\rho\|_{L^2}^{\frac{2}{3}}\|\rho\|_{\dot H^2}
\|\partial_y\rho\|_{\dot H^2}^{\frac{4}{3}}
\notag\\
&\leq C\|\rho\|_{L^2}^{\frac{5}{6}}\|\rho\|_{\dot H^2}^{\frac{1}{2}}
\|\partial_y\rho\|_{\dot H^2}^{\frac{5}{3}}
+C\|\rho\|_{L^2}^{\frac{2}{3}}\|\rho\|_{\dot H^2}
\|\partial_y\rho\|_{\dot H^2}^{\frac{4}{3}}.
\label{eq:IPM-B4}
\end{align}
Young's inequality in \eqref{eq:IPM-B1}--\eqref{eq:IPM-B4}
therefore gives, for every $\delta>0$,
\begin{align}
\sum_{j=1}^4|G_j|
&\leq\delta\|\partial_y\rho\|_{\dot H^2}^2
+C_\delta\bigl(\|\rho\|_{L^2}^2+\|\rho\|_{L^2}^3
+\|\rho\|_{L^2}^5\bigr)\|\rho\|_{\dot H^2}^3.
\label{eq:IPM-B-sum}
\end{align}

We next apply $\partial_y^2$ to the density equation and take the
$L^2$ inner product with $\partial_y^2\rho$. This gives
\begin{align*}
\frac12\frac{d}{dt}\|\partial_y^2\rho\|_{L^2}^2
+\|\partial_y^3\rho\|_{L^2}^2
&=-\int_{\R^2}\partial_y^2u^1\partial_x\rho\partial_y^2\rho\,dx\,dy
-\int_{\R^2}\partial_y^2u^2\partial_y\rho\partial_y^2\rho\,dx\,dy
\\*
&\quad-2\int_{\R^2}\partial_yu^1\partial_{xy}\rho\partial_y^2\rho\,dx\,dy
-2\int_{\R^2}\partial_yu^2(\partial_y^2\rho)^2\,dx\,dy
\\*
&=:H_1+H_2+H_3+H_4.
\end{align*}
Using \eqref{eq:IPM-directional-multipliers} and
\eqref{eq:IPM-interpolation}, we find
\begin{align}
|H_1|
&\leq\|\partial_x\rho\|_{L^2}
\|\partial_y^2u^1\|_{L_y^2L_x^\infty}
\|\partial_y^2\rho\|_{L_x^2L_y^\infty}
\notag\\*
&\leq C\|\partial_x\rho\|_{L^2}
\|\partial_y^2\rho\|_{L^2}^{\frac{1}{2}}
\|\partial_y^3\rho\|_{L^2}^{\frac{1}{2}}
\|\partial_y^2\rho\|_{L^2}^{\frac{1}{2}}
\|\partial_y^3\rho\|_{L^2}^{\frac{1}{2}}
\notag\\*
&\leq C\|\rho\|_{L^2}^{\frac{5}{6}}\|\rho\|_{\dot H^2}^{\frac{1}{2}}
\|\partial_y\rho\|_{\dot H^2}^{\frac{5}{3}}.
\label{eq:IPM-D1}
\end{align}
Similarly, using \eqref{eq:IPM-L2-multiplier},
\begin{align}
|H_2|
&\leq\|\partial_y\rho\|_{L^2}
\|\partial_y^2u^2\|_{L_y^2L_x^\infty}
\|\partial_y^2\rho\|_{L_x^2L_y^\infty}
\notag\\*
&\leq C\|\partial_y\rho\|_{L^2}
\|\partial_y^2\rho\|_{L^2}
\|\partial_{xyy}\rho\|_{L^2}^{\frac{1}{2}}
\|\partial_y^3\rho\|_{L^2}^{\frac{1}{2}}
\notag\\*
&\leq C\|\rho\|_{L^2}\|\partial_y\rho\|_{\dot H^2}^2.
\label{eq:IPM-D2}
\end{align}
For $H_3$, anisotropic H\"older's inequality gives
\begin{align}
|H_3|
&\leq C\|\partial_{xy}\rho\|_{L^2}
\|\partial_yu^1\|_{L_y^2L_x^\infty}
\|\partial_y^2\rho\|_{L_x^2L_y^\infty}
\notag\\*
&\leq C\|\partial_{xy}\rho\|_{L^2}
\|\partial_y\rho\|_{L^2}^{\frac{1}{2}}
\|\partial_y^2\rho\|_{L^2}^{\frac{1}{2}}
\|\partial_y^2\rho\|_{L^2}^{\frac{1}{2}}
\|\partial_y^3\rho\|_{L^2}^{\frac{1}{2}}
\notag\\*
&\leq C\|\rho\|_{L^2}\|\partial_y\rho\|_{\dot H^2}^2.
\label{eq:IPM-D3}
\end{align}
For $H_4$, the same placement of the velocity factor gives
\begin{align}
|H_4|
&\leq C\|\partial_yu^2\|_{L_y^2L_x^\infty}
\|\partial_y^2\rho\|_{L^2}
\|\partial_y^2\rho\|_{L_x^2L_y^\infty}
\notag\\*
&\leq C\|\partial_y\rho\|_{L^2}^{\frac{1}{2}}
\|\partial_{xy}\rho\|_{L^2}^{\frac{1}{2}}
\|\partial_y^2\rho\|_{L^2}^{\frac{3}{2}}
\|\partial_y^3\rho\|_{L^2}^{\frac{1}{2}}
\notag\\*
&\leq C\|\rho\|_{L^2}\|\partial_y\rho\|_{\dot H^2}^2.
\label{eq:IPM-D4}
\end{align}

Apply Young's inequality also to \eqref{eq:IPM-D1}, and choose
$\delta$ sufficiently small in this estimate and
\eqref{eq:IPM-B-sum}. By \eqref{eq:IPM-density-L2},
$\|\rho(t)\|_{L^2}\leq\varepsilon\leq1$.
Taking $\varepsilon_m$ sufficiently small absorbs
\eqref{eq:IPM-D2}--\eqref{eq:IPM-D4} into the vertical
dissipation. Combining the two derivative identities, we obtain
\begin{align}
\frac{d}{dt}\bigl(\|\partial_x^2\rho\|_{L^2}^2
+\|\partial_y^2\rho\|_{L^2}^2\bigr)
+c\|\partial_y\rho\|_{\dot H^2}^2
&\leq C\|\rho\|_{L^2}^2\|\rho\|_{\dot H^2}^3.
\label{eq:IPM-H2-differential}
\end{align}
Here $c>0$ is an absolute constant.

\subsection{Weighted decay and closure of the estimates}

The coefficient on the right-hand side of
\eqref{eq:IPM-H2-differential} will be integrable once we establish
the weighted density decay. We use
\eqref{eq:common-weighted-energy} with $f=\rho$ and $v=u$.
The reflection
\[
(\rho,u^1,u^2)(t,x,y)\longmapsto(-\rho,u^1,-u^2)(t,x,-y)
\]
preserves \eqref{eq:IPM-system}. By uniqueness, $\rho$ and its weighted
counterpart $W$ remain odd in $y$, and
\eqref{eq:moving-coercivity} applies.

We close the energy and decay estimates simultaneously. Fix
$C_0\geq1$, depending only on $m$, to dominate the coefficient of
$\|\rho_0\|_{H^2}^2$ obtained by integrating
\eqref{eq:IPM-H2-differential} and adding
\eqref{eq:IPM-density-L2}. The dependence on this fixed choice is
included in $C_m$ below. For $T>0$ in the interval of existence,
assume
\begin{align}
\sup_{0\leq t\leq T}\|\rho(t)\|_{H^2}^2
+\int_0^T\|\partial_y\rho(t)\|_{H^2}^2\,dt
&\leq4C_0\varepsilon^2,
\label{eq:IPM-energy-bootstrap}\\*
\sup_{0\leq t\leq T}(1+t)^{\frac{3m}{4(m+1)}}\|W(t)\|_{L^2}
&\leq2\|w_m\rho_0\|_{L^2}.
\label{eq:IPM-weighted-bootstrap}
\end{align}
We shall first improve the factor $2$ in the weighted bound to
$3/2$, and then use this decay to improve the energy bound from
$4C_0\varepsilon^2$ to $2C_0\varepsilon^2$.

In the weighted transport term, a uniform $H^2$ bound on the density
alone would leave a nonintegrable coefficient of order
$(1+t)^{-1/2}$. To exploit the vertical regularity supplied by the
dissipation, we use Darcy's law to write
\[
u^2=-\rho+(u^2+\rho).
\]
The first term is the density itself and can be estimated in mixed
norms. The multiplier of the second term is
$\ell^2/(k^2+\ell^2)$, which decays in the horizontal Fourier
variable. This allows an $L^\infty$ bound using only vertical
regularity of the density.

For $L>0$, Fourier inversion and Cauchy--Schwarz give
\begin{align*}
\|u^2+\rho\|_{L^\infty}
&\leq C\int_{\R^2}\frac{\ell^2}{k^2+\ell^2}
|\widehat\rho(k,\ell)|\,dk\,d\ell\\
&\leq C\left(\int_{\R^2}
\frac{\ell^4}{(k^2+\ell^2)^2\left(1+\frac{\ell^4}{L^4}\right)}
\,dk\,d\ell\right)^{\frac12}
\left(\|\rho\|_{L^2}^2+L^{-4}\|\partial_y^2\rho\|_{L^2}^2\right)^{\frac12}.
\end{align*}
The integral in the first factor can be evaluated explicitly:
\begin{align*}
\int_\R\frac{\ell^4}{(k^2+\ell^2)^2}\,dk
&=\frac\pi2|\ell|,\\
\int_{\R^2}
\frac{\ell^4}{(k^2+\ell^2)^2\left(1+\frac{\ell^4}{L^4}\right)}
\,dk\,d\ell
&=\frac\pi2\int_\R\frac{|\ell|}{1+\frac{\ell^4}{L^4}}\,d\ell
=\frac{\pi^2}{4}L^2.
\end{align*}
Consequently,
\[
\|u^2+\rho\|_{L^\infty}
\leq C\bigl(L\|\rho\|_{L^2}
+L^{-1}\|\partial_y^2\rho\|_{L^2}\bigr).
\]
Optimizing in $L$ and then applying
\eqref{eq:IPM-vertical-interpolation}, we obtain
\begin{equation}\label{eq:IPM-vertical-velocity-infinity}
\|u^2+\rho\|_{L^\infty}
\leq C\|\rho\|_{L^2}^{\frac12}\|\partial_y^2\rho\|_{L^2}^{\frac12}
\leq C\|\rho\|_{L^2}^{\frac23}
\|\partial_y\rho\|_{\dot H^2}^{\frac13}.
\end{equation}

\Needspace{9\baselineskip}
We now estimate the two contributions to the weighted transport.
The elementary inequality
\[
2(2m-1)(1+t)+y^2
\geq2\sqrt{2(2m-1)(1+t)}\,|y|
\]
implies
\begin{equation}\label{eq:IPM-weight-coefficient}
\left|\frac{my}{2(2m-1)(1+t)+y^2}\right|
\leq\frac{m}{2\sqrt{2(2m-1)}}(1+t)^{-\frac12}.
\end{equation}
For the contribution containing $-\rho$, anisotropic H\"older's
inequality and interpolation give
\begin{align}
&\left|\int_{\R^2}\frac{my\rho}{2(2m-1)(1+t)+y^2}|W|^2\,dx\,dy\right|\notag\\
\leq{} & C_m(1+t)^{-\frac12}
\|\rho\|_{L_y^2L_x^\infty}\|W\|_{L_x^2L_y^\infty}\|W\|_{L^2}\notag\\
\leq{} & C_m(1+t)^{-\frac12}
\|\rho\|_{L^2}^{\frac12}\|\partial_x\rho\|_{L^2}^{\frac12}
\|W\|_{L^2}^{\frac32}\|\partial_yW\|_{L^2}^{\frac12}\notag\\
\leq{} & C_m(1+t)^{-\frac12}
\|\rho\|_{L^2}^{\frac34}\|\rho\|_{\dot H^2}^{\frac14}
\|W\|_{L^2}^{\frac32}\|\partial_yW\|_{L^2}^{\frac12}.
\label{eq:IPM-weighted-local-intermediate}
\end{align}
By \eqref{eq:moving-coercivity},
\[
\|W\|_{L^2}^{\frac32}
\leq C_m(1+t)^{\frac34}
\left(\|\partial_yW\|_{L^2}^2
+\frac1{1+t}\int_{\R^2}V_m\!\left(\frac{y}{\sqrt{1+t}}\right)
|W|^2\,dx\,dy\right)^{\frac34}.
\]
The derivative factor in
\eqref{eq:IPM-weighted-local-intermediate} is bounded by the
one-quarter power of the same quadratic form. Hence
\begin{align}
&\left|\int_{\R^2}\frac{my\rho}{2(2m-1)(1+t)+y^2}|W|^2\,dx\,dy\right|\notag\\
\leq{} & C_m(1+t)^{\frac14}
\|\rho\|_{L^2}^{\frac34}\|\rho\|_{\dot H^2}^{\frac14}
\left(\|\partial_yW\|_{L^2}^2
+\frac1{1+t}\int_{\R^2}V_m\!\left(\frac{y}{\sqrt{1+t}}\right)
|W|^2\,dx\,dy\right).
\label{eq:IPM-weighted-local}
\end{align}
For the remaining contribution,
\eqref{eq:IPM-vertical-velocity-infinity} and
\eqref{eq:IPM-weight-coefficient} give
\begin{align}
&\left|\int_{\R^2}
\frac{my(u^2+\rho)}{2(2m-1)(1+t)+y^2}|W|^2\,dx\,dy\right|\notag\\
\leq{} & C_m(1+t)^{-\frac12}\|u^2+\rho\|_{L^\infty}\|W\|_{L^2}^2\notag\\
\leq{} & C_m(1+t)^{-\frac12}\|\rho\|_{L^2}^{\frac23}
\|\partial_y\rho\|_{\dot H^2}^{\frac13}\|W\|_{L^2}^2.
\label{eq:IPM-weighted-nonlocal}
\end{align}

We now improve the weighted bound in
\eqref{eq:IPM-weighted-bootstrap}. Since $w_m\geq1$, the two
bootstrap assumptions give
\begin{align}
C_m(1+t)^{\frac14}\|\rho\|_{L^2}^{\frac34}
\|\rho\|_{\dot H^2}^{\frac14}
&\leq C_m\|w_m\rho_0\|_{L^2}^{\frac34}
\varepsilon^{\frac14}
(1+t)^{\frac14-\frac{9m}{16(m+1)}}\notag\\
&\leq C_m\varepsilon
(1+t)^{-\frac{5m-4}{16(m+1)}}.
\label{eq:IPM-small-coefficient}
\end{align}
Choose $\varepsilon_m$ so that this is at most $1/2$.
We can therefore absorb \eqref{eq:IPM-weighted-local} into the
quadratic form in \eqref{eq:common-weighted-energy} and apply
\eqref{eq:moving-coercivity}. Together with
\eqref{eq:IPM-weighted-nonlocal}, this gives
\begin{equation}\begin{aligned}
\frac12\frac{d}{dt}\|W\|_{L^2}^2
+\Biggl[&\frac{3m}{4(m+1)(1+t)}
\left(1-C_m(1+t)^{\frac14}\|\rho\|_{L^2}^{\frac34}
\|\rho\|_{\dot H^2}^{\frac14}\right)\\
&{}-C_m(1+t)^{-\frac12}\|\rho\|_{L^2}^{\frac23}
\|\partial_y\rho\|_{\dot H^2}^{\frac13}
\Biggr]\|W\|_{L^2}^2\leq0.
\end{aligned}\label{eq:IPM-weighted-differential}\end{equation}
The two corrections to $3m/[4(m+1)(1+t)]$ are integrable in time.
For the first, \eqref{eq:IPM-small-coefficient} yields
\begin{align*}
\int_0^T(1+t)^{-\frac34}\|\rho\|_{L^2}^{\frac34}
\|\rho\|_{\dot H^2}^{\frac14}\,dt
&\leq C_m\varepsilon\int_0^\infty
(1+t)^{-1-\frac{5m-4}{16(m+1)}}\,dt\\
&\leq C_m\varepsilon.
\end{align*}
For the second, H\"older's inequality with time exponents $6$ and
$6/5$ uses precisely the integrated dissipation from
\eqref{eq:IPM-energy-bootstrap}:
\begin{align*}
&\int_0^T(1+t)^{-\frac12}\|\rho\|_{L^2}^{\frac23}
\|\partial_y\rho\|_{\dot H^2}^{\frac13}\,dt\\
\leq{} &
\left(\int_0^T\|\partial_y\rho\|_{\dot H^2}^2\,dt\right)^{\frac16}
\left(\int_0^T(1+t)^{-\frac35}\|\rho\|_{L^2}^{\frac45}\,dt\right)^{\frac56}\\
\leq{} & C_m\varepsilon^{\frac13}
\|w_m\rho_0\|_{L^2}^{\frac23}
\left(\int_0^\infty(1+t)^{-\frac35-\frac{3m}{5(m+1)}}\,dt\right)^{\frac56}\\
\leq{} & C_m\varepsilon.
\end{align*}
The last integral is finite because $m>2$. Gronwall's inequality
applied directly to \eqref{eq:IPM-weighted-differential}, together
with these two integral bounds, gives
\begin{equation*}
\|W(t)\|_{L^2}
\leq\|w_m\rho_0\|_{L^2}
(1+t)^{-\frac{3m}{4(m+1)}}\exp(C_m\varepsilon).
\end{equation*}
Taking $\varepsilon_m$ so that the exponential is at most $3/2$
improves \eqref{eq:IPM-weighted-bootstrap}.

This weighted bound now supplies the missing coefficient in the
energy estimate:
\begin{align*}
\int_0^T\|\rho(t)\|_{L^2}^2\,dt
&\leq C_m\|w_m\rho_0\|_{L^2}^2
\int_0^\infty(1+t)^{-\frac{3m}{2(m+1)}}\,dt\\
&\leq C_m\|w_m\rho_0\|_{L^2}^2.
\end{align*}
Integrating \eqref{eq:IPM-H2-differential} and adding
\eqref{eq:IPM-density-L2}, we use
\eqref{eq:IPM-norm-equivalences} to obtain
\begin{align}
&\sup_{0\leq t\leq T}\|\rho(t)\|_{H^2}^2
+\int_0^T\|\partial_y\rho(t)\|_{H^2}^2\,dt\notag\\
\leq{} & C\|\rho_0\|_{H^2}^2
+C\sup_{0\leq t\leq T}\|\rho(t)\|_{\dot H^2}^3
\int_0^T\|\rho(t)\|_{L^2}^2\,dt\notag\\
\leq{} & C\varepsilon^2+C_m\varepsilon^5.
\label{eq:IPM-energy-improvement}
\end{align}
By the choice of $C_0$, the first term is at most
$C_0\varepsilon^2$. Decreasing $\varepsilon_m$ so that
$C_m\varepsilon_m^3\leq C_0$ makes the right-hand side of
\eqref{eq:IPM-energy-improvement} at most
$2C_0\varepsilon^2$, improving
\eqref{eq:IPM-energy-bootstrap}. Both estimates are therefore
closed by continuity, proving \eqref{eq:IPM-main-energy} and
\eqref{eq:IPM-main-decay}. Finally,
\[
\|u(t)\|_{L^2}\leq\|\rho(t)\|_{L^2}
\leq C_m\|w_m\rho_0\|_{L^2}
(1+t)^{-\frac{3m}{4(m+1)}}
\]
by \eqref{eq:IPM-L2-multiplier}, which proves
\eqref{eq:IPM-main-velocity}.

\section*{Acknowledgments}

N. Zhu was partially supported by the National Natural Science Foundation
of China (Grant Nos. 12301285 and 12171010) and the Shandong Provincial
Natural Science Foundation (Grant No. ZR2023QA002).

\vskip .1in
\noindent\textbf{Data availability statement.} Data sharing is not applicable to this article, as no datasets were generated or analyzed during the current study.

\vskip .1in
\noindent\textbf{Conflict of interest.} The authors declare that they have no conflict of interest.


\vskip .3in 

\begin{thebibliography}{99}
\small
\setlength{\itemsep}{1pt}

\bibitem{CG}
C. Cao and Y. Guo,
On the two-dimensional Navier--Stokes equations with horizontal viscosity,
\textit{J. Math. Anal. Appl.} \textbf{557} (2026), Paper No. 130352.

\bibitem{CDGG}
J.-Y. Chemin, B. Desjardins, I. Gallagher and E. Grenier,
Fluids with anisotropic viscosity,
\textit{ESAIM Math. Model. Numer. Anal.} \textbf{34} (2000), 315--335.

\bibitem{CZ}
J.-Y. Chemin and P. Zhang,
On the global wellposedness to the 3-D incompressible anisotropic Navier--Stokes equations,
\textit{Comm. Math. Phys.} \textbf{272} (2007), 529--566.

\bibitem{CJ}
K. Choi and I.-J. Jeong,
Infinite growth in vorticity gradient of compactly supported planar vorticity near Lamb dipole,
\textit{Nonlinear Anal. Real World Appl.} \textbf{65} (2022), Paper No. 103470.

\bibitem{CMIPM}
D. C\'ordoba and L. Mart\'inez-Zoroa,
Finite time singularities of smooth solutions for the 2D incompressible porous media (IPM) equation with a smooth source,
arXiv:2410.22920 (2024).

\bibitem{DaiIPM}
M. Dai,
Finite-time blowup for 2D unforced IPM,
arXiv:2610.06712 (2026).

\bibitem{Dong}
B. Dong, J. Wu, X. Xu and N. Zhu,
Stability and exponential decay for the 2D anisotropic Navier--Stokes equations with horizontal dissipation,
\textit{J. Math. Fluid Mech.} \textbf{23} (2021), Paper No. 100.

\bibitem{Fujii}
M. Fujii,
Large time behavior of solutions to the 3D anisotropic Navier--Stokes equation,
\textit{Nonlinear Anal. Real World Appl.} \textbf{71} (2023), Paper No. 103821.

\bibitem{GW}
T. Gallay and C. E. Wayne,
Invariant manifolds and the long-time asymptotics of the Navier--Stokes and vorticity equations on $\R^2$,
\textit{Arch. Ration. Mech. Anal.} \textbf{163} (2002), 209--258.

\bibitem{JYaoZhou}
I.-J. Jeong, Y. Yao and T. Zhou,
Superlinear gradient growth for 2D Euler equation without boundary,
arXiv:2507.15739 (2025).

\bibitem{JWY}
R. Ji, J. Wu and W. Yang,
Stability and optimal decay for the 3D Navier--Stokes equations with horizontal dissipation,
\textit{J. Differential Equations} \textbf{290} (2021), 57--77.

\bibitem{KS}
A. Kiselev and V. \v{S}ver\'ak,
Small scale creation for solutions of the incompressible two-dimensional Euler equation,
\textit{Ann. of Math.} \textbf{180} (2014), 1205--1220.

\bibitem{KY}
A. Kiselev and Y. Yao,
Small scale formations in the incompressible porous media equation,
\textit{Arch. Ration. Mech. Anal.} \textbf{247} (2023), Paper No. 1.

\bibitem{LZZ}
S. Liang, P. Zhang and R. Zhu,
Deterministic and stochastic 2D Navier--Stokes equations with anisotropic viscosity,
\textit{J. Differential Equations} \textbf{275} (2021), 473--508.

\bibitem{Paicu}
M. Paicu,
\'Equation anisotrope de Navier--Stokes dans des espaces critiques,
\textit{Rev. Mat. Iberoam.} \textbf{21} (2005), 179--235.

\bibitem{PZ}
M. Paicu and P. Zhang,
Global solutions to the 3-D incompressible anisotropic Navier--Stokes system in the critical spaces,
\textit{Comm. Math. Phys.} \textbf{307} (2011), 713--759.

\bibitem{Pran}
L. Prandtl,
Zur Berechnung der Grenzschichten,
\textit{Z. Angew. Math. Mech.} \textbf{18} (1938), 77--82.

\bibitem{Sch}
H. Schlichting and K. Gersten,
\textit{Boundary-layer theory}, ninth edition,
translated from the German by Katherine Mayes, Springer, Berlin, 2017.

\bibitem{Schonbek}
M. E. Schonbek,
$L^2$ decay for weak solutions of the Navier--Stokes equations,
\textit{Arch. Rational Mech. Anal.} \textbf{88} (1985), 209--222.

\bibitem{SCho4}
M. E. Schonbek,
Large time behaviour of solutions to the Navier--Stokes equations in $H^m$ spaces,
\textit{Comm. Partial Differential Equations} \textbf{20} (1995), 103--117.

\bibitem{SW}
M. E. Schonbek and M. Wiegner,
On the decay of higher-order norms of the solutions of Navier--Stokes equations,
\textit{Proc. Roy. Soc. Edinburgh Sect. A} \textbf{126} (1996), 677--685.

\bibitem{WangWuZhu}
Z. Wang, J. Wu and N. Zhu,
Stability and decay for the 2D anisotropic Navier--Stokes equations with fractional horizontal dissipation on $\R^2$,
arXiv:2601.16110 (2026).

\bibitem{Wiegner}
M. Wiegner,
Decay results for weak solutions of the Navier--Stokes equations on $\R^n$,
\textit{J. London Math. Soc.} (2) \textbf{35} (1987), 303--313.

\bibitem{XZ}
L. Xu and P. Zhang,
Enhanced dissipation for the third component of 3D anisotropic Navier--Stokes equations,
\textit{J. Differential Equations} \textbf{335} (2022), 464--496.

\bibitem{ZhengIPM}
F. Zheng,
Finite-time blowup of the IPM equation from smooth initial data without forcing,
arXiv:2610.06544 (2026).

\bibitem{ZW}
D. Zhou and J. Wu,
$H^1$-uniqueness, stability and decay for anisotropic Navier--Stokes and Boussinesq equations,
\textit{J. Differential Equations} \textbf{453} (2026), Paper No. 113930.

\bibitem{ZlatosHalfPlane}
A. Zlato\v{s},
Maximal double-exponential growth for the Euler equation on the half-plane,
\textit{Invent. Math.} \textbf{243} (2026), 117--126.

\end{thebibliography}
\end{document}